\documentclass[11pt]{article}
\usepackage[letterpaper,margin=1in]{geometry}
\usepackage{amsmath,amssymb,amsthm,mathtools}
\usepackage{booktabs,graphicx,array}
\usepackage{microtype}
\usepackage{natbib}
\usepackage{soul}
\usepackage{algorithm}
\usepackage{algpseudocode}
\usepackage{tikz}
\usepackage[most]{tcolorbox}
\usepackage[hidelinks]{hyperref}
\newtheorem{theorem}{Theorem}
\theoremstyle{remark}
\newtheorem{remark}{Remark}
\newcommand{\native}{P_0}
\newcommand{\X}{\mathcal{X}}
\newcommand{\mX}{\mathcal{X}}
\newcommand{\mP}{\mathcal{P}}

\makeatletter
\newcommand{\namedlabel}[2]{%
  \begingroup
  \def\@currentlabel{#2}%
  \label{#1}%
  \endgroup
}
\makeatother

\newcommand{\var}{\mbox{Var}}

\newcommand{\E}{\mathbb{E}}
\newcommand{\mc}{\hat{\theta}_{\mathrm{MC}}}
\newcommand{\cv}{\hat{\theta}_{\mathrm{cv}}}
\newcommand{\lr}{\hat{\theta}_{\mathrm{lr}}}
\newcommand{\cvbudget}{\bar{\theta}_{\mathrm{cv}}}
\newcommand{\lrbudget}{\bar{\theta}_{\mathrm{lr}}}

\newcommand{\trans}{\top} 
\newcommand{\simout}{Y}
\newcommand{\simoutXQ}{Y(X,Q)}

\newcommand{\simoutdotQ}{Y(\cdot,Q)}
\newcommand{\simoutQ}{Y_Q}

\newcommand{\eat}[1]{}

\begin{document}
\title{Combining Control Variates and Likelihood Ratios for Variance Reduction in Simulation}
\author{Shane G. Henderson\\
\small School of Operations Research and Information Engineering, Cornell University
\and Raghu Pasupathy\\
\small Department of Statistics, Purdue University}
\date{}
\maketitle
\begin{abstract}
When using (stochastic) simulation to estimate an expected value of some random quantity, one must specify both the random quantity that is the output of the simulation logic and an underlying probability distribution. The variance-reduction method of control variates modifies the simulation output to obtain a lower-variance estimator with the same mean, while the method of importance sampling modifies the underlying probability distribution. We develop a method for combining these approaches. An initial optimization process alternates between choosing the coefficients for the control variates and updating the probability distribution iteratively to attempt to minimize the variance of a combined estimator. We explore both nonparametric and parametric specifications and updates of the probability distribution. It is possible to implement the new estimators while still sampling from the initial probability distribution through judicious use of likelihood ratios, which can ease implementation. Estimation then follows in a separate step using the optimized combined estimator. Examples demonstrate that the combined approach can yield highly efficient estimators, though sometimes the gains are modest relative to the use of just the better of control variates and importance sampling.
\end{abstract}
\par\smallskip
\noindent\textbf{Abbreviations:} MC, Monte Carlo; CV, control variates; LR, likelihood ratio; CVLR, control variates likelihood ratio; AU-CVLR, alternating updates CVLR.
\par\smallskip
\noindent\textbf{Corresponding author:} Raghu Pasupathy (\href{mailto:pasupath@purdue.edu}{pasupath@purdue.edu}).

\section{Introduction}\label{sec1}

Consider the classical simulation problem of estimating the expectation
\begin{equation}\label{theta}
\theta  = \mathbb{E}_{X \sim \native}[\simout_0(X)] = \int_{x \in \mathcal{X}} \simout_0(x)\,P_0(dx), 
\end{equation}
where $P_0\in\mathcal P(\X)$ is the ``native'' probability distribution describing the random quantities underlying the simulation, $\simout_0:\X\to\mathbb R$ is the corresponding ``native'' simulation output function, and $\mathcal{X}$ is the support. As a concrete example, $\theta$ might denote the expected waiting time of the twentieth customer in an M/M/1 queue that starts idle; $\simout_0(x)$ in~\eqref{theta} then corresponds to the waiting time of the twentieth customer in the specific ``scenario'' $x \in \mathbb{R}^{2 \times 20}$ that specifies the inter-arrival and service times of the first 20 customers, and $P_0$ the sampling law governing the generation of the scenario $x \in \mathcal{X}$. 

The crude Monte Carlo (CMC) estimator of $\theta$ is the classical ``one-sample'' estimator constructed using the native simulation, \begin{equation}\label{cmc}\hat{\theta}_{\mathrm{MC}} = \simout_0(X), \quad X \sim \native \tag{CMC}\end{equation} 
that is then replicated, independently, many times and averaged to estimate $\theta$. The estimator $\mc$ is an unbiased estimator of $\theta$ because $\mathbb{E}[\mc] = \int \simout_0(x) P_0(dx) = \theta$, and has variance \begin{equation}\label{varmc} \mbox{Var}(\mc) = \int \simout_0^2(x) P_0(dx) - \theta^2.
\end{equation} 

Now suppose we have the ability to alter the native simulation $\simout_0(\cdot)$ and the native law $P_0$. Then we might ask: can the ability to alter $\simout_0$ and $P_0$ be used toward \emph{variance reduction}, that is, constructing estimators of $\theta$ that have variance smaller than $\mbox{Var}(\mc)$? 

Variance reduction is often possible when exogeneous information is present, and such information can be used to modify the estimator of $\theta$. For instance, suppose we have external information available in the form of a known (population) mean vector $\gamma = \int g(x) \native (dx)$ with respect to a known ($\native$-measurable) function $g: \mathcal{X} \to \mathbb{R}^d$. In the M/M/1 example, knowing the vector giving the mean interarrival time and mean service time, $\gamma = (\lambda^{-1}, \mu^{-1})^\trans$ (with $\trans$ denoting the transpose operation), we might take $g$ as the function that returns the (sample) average of the first 20 interarrival and service times. Then a classical way of incorporating such information is through the \emph{control variates} (CV) estimator,  \begin{equation}\label{cv} \hat{\theta}_{\mathrm cv} = \simout_0(X) - \beta^\top(g(X) - \gamma), \quad X \sim \native, \tag{CV}\end{equation}
where the parameter $\beta$ is analogous to the slope parameter in classical linear regression and can be chosen through a variety of means, primarily estimation. The estimator $\cv$ has been studied and surveyed extensively, e.g.,~\cite{nel90,glysze02} and \cite[\S4.1]{Glasserman04}; for instance, it is known to be an unbiased estimator of $\theta$ with lower variance than $\mc$ under very general conditions.

A second and key example for this paper arises when we have the option of simulating under an alternative probability distribution, $Q$. If we do that without modifying the simulation output $Y_0$ then bias will typically result, in that
$$\E_{X \sim Q}[Y_0(X)] = \int_{x\in \mathcal X} Y_0(x) Q(dx) \ne \theta.$$
We can recover an unbiased estimator by multiplying by the {\em likelihood ratio}, $\frac{dP_0}{dQ}(x)$, which is well defined under the condition that $Q \gg \native$. The notation $Q \gg \native$ is the standard ``absolute continuity'' condition that avoids the situation where $Q$ assigns zero probability to an event that has positive probability under $\native$. The ratio $\frac{dP_0}{dQ}(x)$ is the Radon-Nikodym derivative, which is the ratio of probabilities when both $P_0$ and $Q$ are discrete probability distributions and the ratio of densities when $P_0$ and $Q$ are continuous distribution functions with densities. The likelihood ratio estimator, which is unbiased, is then
\begin{equation}\label{lr} \hat{\theta}_{\mathrm lr} = \simout_0(X) \frac{d\native}{dQ}(X), \quad X \sim Q, \tag{LR}
\end{equation}
and has also been studied extensively; see, e.g., Section V.1 and Chapter~VI of \cite{asmgly07}.

The simulation estimators $\mc, \cv$ and $\lr$ are all constructed using (only) a pair of mathematical objects:
\begin{enumerate}
\item[(i)] a \emph{sampling law} that we denote $Q \in \mathcal{P}(\mathcal{X})$, where $\mathcal P(\mX)$ is the set of probability distributions on $\mX$, and \item[(ii)] a \emph{simulation output function} that we denote $Y(\cdot,Q): \mathcal{X} \to \mathbb{R}.$
\end{enumerate}
Then, making choices for $Q$ and $\simoutQ = \simoutdotQ$ results in the one-sample estimator $$\simoutXQ, \quad X \sim Q$$ for the quantity $\theta$. Thus, the CMC estimator $\mc$ in~\eqref{cmc} chooses $(Q,\simoutQ):= (P_0,\simout_0)$, the CV estimator $\cv$ in~\eqref{cv} chooses $(Q,\simoutQ):= (P_0,\simout_0(\cdot) - \beta^\top(g(\cdot) - \gamma))$, and the LR estimator in~\eqref{lr} chooses $(Q, \simoutQ) = (Q, Y_0(\cdot) \frac{dP_0}{dQ}(\cdot))$. We sometimes refer to the pair $(Q,\simoutQ)$ as an \emph{estimator pair}. 

In this paper, we study the {\em CVLR} class of estimator pairs that combines the CV and LR approaches, i.e., that sets
\begin{equation}\label{oracularmap}
\simout_{\beta}(X,Q) := \left(\simout_0(X) - \beta^\top(g(X) - \gamma)\right)\, \underbrace{\frac{d\native}{dQ}(X)}_{L(X)},
\end{equation}
where $X \sim Q$ and $Q \gg \native$. The random variable $\simout_{\beta}(X,Q)$ is then unbiased for any fixed $\beta$ since \begin{align}\label{unbiasedT}
\int \simout_{\beta}(x,Q) Q(dx) &:= \int \left(\simout_0(x) - \beta^\top(g(x) - \gamma)\right)\,L(x) Q(dx) \nonumber \\
&= \int \simout_0(x)\native(dx) - \underbrace{\int (g(x) - \gamma)\native(dx)}_{=0} = \theta.
\end{align}
Minimizing the variance, $\mbox{Var}_{X \sim Q}(\simout_{\beta}(X,Q))$, of the unbiased estimator $\simout_{\beta}(X,Q), X \sim Q$ of $\theta$ thus entails optimizing over the set $\{(Q, \beta) \in \mP(\mX) \times \mathbb{R}^d: Q \gg \native\}$. (The optimal $Q$ is not in general a trivial point-mass distribution on some $x$ for which $\simout_\beta(x, Q) = \theta$ because we have stipulated that $Q \gg \native.$) %

To see why it can be beneficial to combine control variates and importance sampling, consider the following, highly stylized, example. Suppose $Y_0(x) = e^{-x} + \eta I(x \le 1)$ for $x \ge 0$,
where $I(\cdot)$ is the indicator function that equals 1 when its argument is true and 0 otherwise, and $\eta > 0$ is a deterministic constant. Moreover, suppose that $\native$ is the exponential distribution with mean 1. Then $\theta = 1/2 + \eta(1-e^{-1})$.
Assume a single control $g(x) = I(x \le 1)$ with mean
$\gamma = 1-e^{-1}$ and the family of importance sampling distributions given by a mixture of two exponential distributions
$$\mathcal{Q} = \{Q\sim a\exp(r_1) + (1-a) \exp(r_2): r_1, r_2 > 0, a \in [0, 1]\}.$$

When $\eta > 0$, using control variates alone cannot yield 0 variance, nor can the use of importance sampling alone. However, if we take the control variate coefficient $\beta = \eta$ then $Y_\beta(x) = e^{-x} + \eta(1-e^{-1})$
for $x \ge 0$. We can then write
$$\theta = \int_0^\infty Y_\beta(x)e^{-x}\, dx = \int_0^\infty \left[e^{-2x} + \eta(1-e^{-1}) e^{-x}\right] \, dx.$$

Thus, if we take the importance sampling density
$$q(x) = \frac{1}{2\theta} 2e^{-2x} + \frac{2\eta(1-e^{-1})}{2\theta} e^{-x},$$
which lies within the family $\mathcal{Q}$, then we attain 0 variance when sampling under $q$.

In this example, control variates are used to reshape the integrand $Y_0$ to $Y_\beta$ such that $Y_\beta p_0$ (where $p_0$ is the density of $P_0$) is an exact multiple of one of the importance sampling distributions in our class. In general, control variates can be used to reshape the integrand to ensure this holds as closely as possible.

Returning to the general setting, as we shall see, optimizing over $\beta$ is straightforward, reducing to a least squares problem just as it does with standard CV estimators. However, optimization over the probability distribution $Q$ is more challenging. 

One approach, taken in the stylized example above, is to assume a parametric form for the distribution $Q$, then use some custom method to optimize over those parameters. That is the approach taken in an option pricing example detailed in Section~\ref{sec:barrier-asian}, where we use a pilot sample to estimate the variance associated with {\em every} set of parameters from within a finite set and select the set of parameters that gives the smallest estimated variance. That effort could alternatively be performed using ranking and selection techniques \citep{nelson2022let}, but we did not do that due to the added complexity. One could also perform a gradient-based search over the parameter space, which would require estimating gradients. We did not take this approach either.

A second approach, developed in Sections~\ref{sec:CVLRAU} and \ref{sec:semaucvlr}, is to adopt a nonparametric perspective in which we attempt to optimize over a very broad range of distributions $Q$. In that section, we develop an iterative approach where, in each iteration, we optimize $\beta$ for the current distribution $Q_k$, then update $Q_k$ to $Q_{k+1}$ with the newly optimized $\beta$. It is not generally tractable to identify the optimal $Q_{k+1}$ because the variance expression is nonlinear, so we instead {\em linearize} the variance expression using a Taylor expansion. We also employ {\em regularization} to ensure that the distribution $Q_k$ does not change so much that the Taylor expansion becomes too inaccurate. We then show that the (original) joint optimization problem over $\beta$ and $Q$ is a convex optimization problem and establish convergence theory. The distributions $(Q_k: k \ge 1)$ so produced are not trivial to sample from, but we provide both a rejection sampler and a Markov Chain Monte Carlo sampler.

The contributions of this article are as follows.
\begin{enumerate}
    \item We develop an approach to combining control variates and likelihood ratios to achieve variance reduction, showing that these methods can cooperatively reduce variance.
    \item We develop a tractable approach to optimization over a nonparametric family of distributions $Q$. We prove that the resulting (in general infinite-dimensional) optimization problem is convex; see the proof of Theorem~\ref{thm:aemm-optimality}.
    \item Through linearization and regularization we develop an iterative approach to nonparametric optimization and prove that the iterates converge to a global minimizer of the regularized problem.
    \item We develop an implementation of the nonparametric optimization and sampling algorithm and demonstrate it on 2 examples, a somewhat contrived M/M/1 example, and a realistic emergency services example.
    \item We further demonstrate the parametric approach on an option-pricing example, demonstrating important efficiency gains.
\end{enumerate}

The remainder of this paper is organized as follows. Section~\ref{sec:lit} reviews the relevant literature. Section~\ref{sec:CVLRAU} begins our exploration of the nonparametric approach together with its associated theory. Section~\ref{sec:semaucvlr} develops a specific implementation based on mirror descent. Section~\ref{sec:numerical} uses the M/M/1 example to demonstrate the nonparametric approach and the effect of regularization. We also demonstrate the nonparametric approach on an emergency services example, and conclude with the parametric approach on an option pricing example. As we shall see, the nonparametric approach does not always yield strong variance reduction beyond that of the better of CV or LR alone, but substantial efficiency improvement is possible as seen in the option-pricing example.

\section{Literature}
\label{sec:lit}
Importance sampling and control variates are among the most widely used and most general variance-reduction methods in stochastic simulation. For general treatments see Chapters~V and VI of \cite{asmgly07}, or Chapter~4 of \cite{Glasserman04}.

The control-variate method has a well-known least-squares and Hilbert-space interpretation involving projection of the centered $Y_0$ onto the linear span of the centered controls. residual after projection. This geometric view, emphasized by \cite{sze06} as part of a Hilbert space view of many variance reduction methods (see also \cite{glysze02}), connects control variates with approximation theory and suggests the use of increasingly rich collections of basis functions as in \cite{Maire2003}. In particular, the zero-variance control variate can be viewed as an ideal approximation target, motivating basis-function, regression, and adaptive-control-variate methods. Highly related developments include practical remedies and coefficient-estimation issues \citep{nel90}, martingale and approximation-based variance reduction \citep{hengly02}, adaptive control variates \citep{kimhen06}, and more recent regression-based, Stein-based, and control-functional approaches to constructing flexible control variates \citep{mira2013zero,south2023regularized,oates2017control}.

Importance sampling likewise has a zero-variance benchmark. Much of the adaptive importance-sampling literature can therefore be interpreted as approximating this ideal distribution within a parametric, mixture, exponential-family, or sequentially updated class. Important approaches include early work in physics \citep{booth2001approximate}, cross-entropy methods \citep{rubinstein2004cross}, adaptive mixture importance sampling, \citep[for example]{cappe2008adaptive,cornuet2012adaptive,owen2020square}, and defensive and mixture importance sampling \citep{hesterberg1995weighted,owen2000safe,he2014optimal}. The latter paper also showed joint convexity of the variance as a function of the control variate coefficients and mixing weights in a mixed importance sampling setting.

Several works explicitly combine, or closely connect, importance sampling and control variates. \cite{owen2000safe}, for example, show how proposal distributions and defensive mixtures can generate useful control variates and provide protection against unstable importance weights. Hesterberg’s work on defensive importance sampling and subsequent work on mixture importance sampling likewise highlights that proposal selection, weighting, and auxiliary information should be considered jointly rather than as wholly separate design choices.

Related optimization problems arise in multiple importance sampling, where samples from several proposal distributions are combined and the estimator weights or sample allocations are selected to reduce variance. \cite{veach1995optimally} introduced important weighting rules for combining sampling techniques, while later work has studied optimal mixture weights, generalized multiple importance sampling, and optimal allocation across proposals; see, for example, \cite{elvira2019generalized}. Although these methods do not always take the form of classical control variates, they share the central feature of optimizing an estimator over a structured family of unbiased or approximately unbiased variance-reduction schemes.

From this perspective, jointly selecting an importance-sampling distribution and a control-variate coefficient vector is a natural extension of existing adaptive methods.

Although prior work has combined importance sampling with control variates and has developed adaptive procedures for each method separately, relatively little work formulates the simultaneous selection of a control-variate coefficient vector and an importance-sampling distribution within a common finite-dimensional optimization problem. Our formulation makes explicit the dependence of the optimal control coefficients on the proposal distribution and, correspondingly, the dependence of proposal quality on the residual variance after control-variate adjustment.

\section{CVLR and Alternating Updates}
\label{sec:CVLRAU}
We say that an estimator is a CVLR estimator if it is of the form $\simout_\beta(X, Q), X \sim Q$ with $Q \gg \native$. CVLR estimators can be constructed in many ways. We now develop, in detail, a specific approach to constructing CVLR estimators based on solving, exactly or approximately, the regularized optimization problem
\begin{align}\label{vrmaster}
\min_{\beta \in \mathbb{R}^d,Q} :\quad & \int \simout_{\beta}(x,Q)^2 \, Q(dx) + \lambda D_{\mathrm{KL}}(Q \parallel \native )\nonumber \\
\mathrm{subject \,\, to} : \quad & Q \equiv \native, \quad Q \in \mathcal{P}(\mathcal{X}) \tag{VR-master}
\end{align}
where the notation $Q \equiv \native$ means the two distributions are {\em equivalent} in the sense that $Q \gg \native$ and $\native \gg Q$, \(\lambda\ge 0\) is a regularization parameter, $$D_{\mathrm{KL}}(\nu \parallel \mu) := \int \log \frac{d\nu}{d\mu}(x) \, d\nu(x), \nu \ll \mu $$
is the KL-divergence between the probability measures $\nu$ and $\mu$, and $d\nu/d\mu$ is the Radon-Nikodym derivative of $\nu$ with respect to $\mu$. We have already observed that the estimator $\simout_{\beta}(X,Q), X \sim Q$ is unbiased with respect to $\theta$. Therefore, minimizing its variance is equivalent to minimizing the second moment, which is the first term in the objective function of~\eqref{vrmaster}. The second term is a regularization term~\cite{2017bec} that penalizes departures of $Q$ from the native law $\native$. This term requires the additional condition that $\native \gg Q$, which we therefore impose. Regularization terms are often used in optimization to increase numerical stability. Thus, if $(\beta_k,Q_k)$ is the solution obtained after executing $k$ iterations of an algorithm to solve~\eqref{vrmaster}, the corresponding CVLR estimator of $\theta$ is $\simout_{\beta_k}(X,Q_k), X \sim Q_k$.

\subsection{CVLR versus Pure-CV, Pure-LR}

Some further observations about~\eqref{vrmaster} are salient. First, since $Q=P_0$ trivially satisfies the constraints $Q \equiv P_0, Q \in \mathcal{P}(\mathcal{X})$, the problem in~\eqref{vrmaster} has a non-empty feasible region for any $\lambda \geq 0$. Second, suppose $(\beta_{\lambda}^\ast,Q_{\lambda}^\ast)$ solves~\eqref{vrmaster}. Then, for any $\lambda \geq 0$, the estimator \begin{equation}\label{cvlr}\simout_{\beta_{\lambda}^\ast}(X,Q_{\lambda}^\ast), \quad X \sim Q_{\lambda}^\ast \tag{CVLR}\end{equation} has lower variance than both $\mc$ and the pure CV estimator $\cv= \simout_{{\beta}^\ast_{\lambda=0}}(X,\native), X \sim \native$. Indeed, since
\begin{align*}\left\{(\beta,Q) \in \mathbb{R}^d \times \mathcal{P}(\mathcal{X}): Q \equiv \native\right\} & \supset \left\{(\beta,\native) \in \mathbb{R}^d \times \mathcal{P}(\mathcal{X}): Q \equiv \native\right\} \\ &\supset \left\{(0,\native): Q \equiv \native\right\},
\end{align*}
we have the corresponding ordering
\begin{align*}
\lefteqn{\int \simout_{\beta_{\lambda}^\ast}(x,Q_{\lambda}^\ast)^2 Q_{\lambda}^\ast(dx) + \lambda \underbrace{D_{\mathrm{KL}}(Q_{\lambda}^\ast \parallel \native )}_{\geq 0}} \qquad\qquad \qquad\\
& \leq  \int (\simout_0(x) - \beta^\ast_{\lambda=0}(g(x)-\gamma))^2 \native(dx) \\
& \leq  \int Y^2_0(x) \native(dx) + \lambda \underbrace{D_{\mathrm{KL}}(\native \parallel \native )}_{0}.
\end{align*}
Thus, for any $\lambda \geq 0$, we have the hierarchy \begin{equation}\label{purecvcomp}
\var_{X \sim Q_{\lambda}^\ast}(\simout_{\beta_\lambda}(X,Q_{\lambda}^\ast)) \leq \var(\cv) \leq \var(\mc).
\end{equation}
The ordering in~\eqref{purecvcomp} establishes that solving~\eqref{vrmaster} (for any fixed $\lambda \geq 0)$ results in an estimator $\simout_{\beta_\lambda}(X,Q_{\lambda}^\ast), X \sim Q^*_{\lambda}$ that dominates the pure-CV estimator $\cv$, which in turn dominates the crude MC estimator $\mc.$ 

Third, the regularization parameter $\lambda$ also serves as a formal mechanism for combining pure CV and pure importance sampling. Specifically, $\lambda \to 0$ promotes aggressive redesign of the sampling law $Q$ in pursuit of variance reduction, whereas, $\lambda \to \infty$ promotes the redesign of $\simoutQ$ while ensuring that $Q$ stays close to the native law $\native$.

Observe that the pure-LR estimator does not appear in the hierarchy~\eqref{purecvcomp}. While pure-LR is easily seen to have variance no larger than that of the CMC estimator, a clean comparison with the CVLR and pure-CV estimators is obtained through the following KL-budget reformulation of~\eqref{vrmaster}:
\begin{align}
\label{vrmasterreform}
\min_{\beta \in \mathbb{R}^d,Q} \quad
& \int \simout_{\beta}(x,Q)^2 \, Q(dx)
\nonumber\\
\mathrm{subject\ to}\quad
& D_{\mathrm{KL}}(Q \parallel \native ) \leq \rho,
\nonumber\\
& Q \sim \native, \qquad Q \in \mathcal{P}(\mathcal{X}).
\tag{KL-budget}
\end{align}
Here \(\rho \geq 0\) is a prescribed KL budget.

The formulations~\eqref{vrmaster} and~\eqref{vrmasterreform} describe the same efficient frontier between second moment and KL displacement. In particular, for each fixed \(\lambda \geq 0\), if
\((\beta^\star_{\lambda},Q^\star_{\lambda})\) solves~\eqref{vrmaster}, then it also solves~\eqref{vrmasterreform} with
\[
\rho_{\lambda}
=
D_{\mathrm{KL}}(Q^\star_{\lambda}\parallel\native).
\]
Conversely, under the usual convexity and constraint-qualification conditions, a solution
\((\beta^\star_{\rho},Q^\star_{\rho})\) of~\eqref{vrmasterreform} is also a solution of~\eqref{vrmaster} for an associated Lagrange multiplier
\(\lambda_{\rho}\geq0\). Thus \(\lambda\) may be viewed as the price of KL displacement, while \(\rho\) specifies the available KL budget.

The formulation~\eqref{vrmasterreform} lends itself to a direct comparison of the variances achieved by the four estimators. Specifically, suppose
\((\bar\beta_{\rho},\bar Q_{\rho})\) solves~\eqref{vrmasterreform}, and let
\(\lrbudget\) and \(\cvbudget\) denote the corresponding solutions obtained by imposing, respectively, \(\beta=0\) and \(\rho=0\). The first restriction yields pure LR, whereas the second forces \(Q=\native\) and hence yields pure CV. Since both restricted problems have feasible sets contained in that of~\eqref{vrmasterreform}, and since CMC is itself feasible for both, simple nesting arguments give
\begin{equation}
\label{allcomp}
\var_{X \sim \bar Q_{\rho}}
\!\left(
\simout_{\bar \beta_\rho}(X,\bar Q_{\rho})
\right)
\leq
\min\!\left\{
\var(\cvbudget),
\var(\lrbudget)
\right\}
\leq
\var(\mc).
\end{equation}

\subsection{Alternating Updates CVLR}
\label{sec:aemm-md}

We now propose an iterative algorithm called \emph{alternating updates CVLR} (AU-CVLR) that solves \eqref{vrmaster} by performing repeated alternating $\beta \in \mathbb{R}^d$ and $Q \in \mathcal{P}(\mathcal{X})$ updates. Specifically, the AU-CVLR updates have the alternating form

\begin{tcolorbox}[
    enhanced,
    colback=white,
    colframe=black,
    boxrule=0.6pt,
    arc=0pt,
    outer arc=0pt,
    attach boxed title to top left={
        xshift=3mm,
        yshift=-2mm
    },
    boxed title style={
        colback=white,
        colframe=white,
        boxrule=0pt
    },
    title={AU-CVLR},
    fonttitle=\bfseries,
    coltitle=black
]
\namedlabel{eq:aem-measure}{AU-CVLR}
\begin{align}
\beta_{k+1}
& \in
\mathop{\arg\min}_{\beta \in \mathbb{R}^d}
\int \simout_{\beta}(x,Q_k)^2\,Q_k(dx), \nonumber \\
&\qquad \simout_{k+1} := \simout_{\beta_{k+1}}(Q_k) \tag{CV-side} \\
Q_{k+1}
& \in
\mathop{\arg\min}_{Q\in\mathcal Q}
\left\{
\int \simout_{k+1}(x,Q)^2\,Q(dx)
+
\lambda D_{\mathrm{KL}}(Q \parallel \native )
\right\}.   \tag{LR-side}
\end{align}
\end{tcolorbox}

The (implicit) CV-side feasible region $\mathcal{Y}_{\mathrm CV}(Q_k)$ is a $\beta$-parameterized convex set of functions given by
\begin{equation}\label{Qkplusone}
\mathcal{Y}_{\mathrm CV}(Q_k) := \left\{\simout_{\beta}(\cdot,Q_k):= \left(\simout_0(\cdot) - \beta^\top(g(\cdot) - \gamma)\right)\, \frac{d\native}{dQ_k}(\cdot), \,\beta \in \mathbb{R}^d \right\},  \end{equation} whereas the explicit LR-side feasible region $\mathcal{Q}$ is the convex set of probability measures supported on $\mathcal{X}$ and equivalent to $\native$, that is,
\begin{equation} \mathcal Q := \left\{Q \in \mathcal{P}(\mathcal{X}): Q \equiv \native\right\}.
\end{equation} 

The expression for $Y_{\beta}$ in $\mathcal{Y}_{\mathrm CV}$ suggests that we need only $Q \gg \native$. And yet, we assume that $Q \in \mathcal{Q}$ satisfies both $Q \gg \native$ and $\native \gg Q.$ The latter is needed for two reasons. First, we already noted that the regularization term appearing in~\eqref{Qkplusone} involves a KL-divergence term with respect to $\native$. Second, the LR-side update calculations that we present later rely on both $\frac{dQ}{d\native}(x)$ and $\frac{d\native}{dQ}(x)$, with the former dictating where to move mass, thereby restricting the search space.

\begin{remark}[Nonparametric and parametric CVLR updates] The CV-side update in~\eqref{eq:aem-measure} is a \emph{parametric} update involving $\beta \in \mathbb{R}^d$, whereas the LR-side update is \emph{nonparametric} as written, that is, no ``structural'' assumptions are made on the feasible region $\mathcal{Q}$. An algorithm implementing the LR-side update will accordingly be nonparametric, as we describe in Section~\ref{sec:lrside}. By contrast, the feasible region $\mathcal{Q}$ can also be parameterized --- see the options pricing example discussed in Section~\ref{sec:barrier-asian}. An appropriate parametric LR-side update algorithm will then have to be devised.  
\end{remark}

That the feasible regions $\mathcal{Y}_{\mathrm CV}(Q_k)$ and $\mathcal{Q}$ are convex is easy to see since $\simout_{\beta_1}, \simout_{\beta_2} \in \mathcal{Y}_{\mathrm CV}(Q_k) \implies \alpha \simout_{\beta_1} + (1-\alpha)\simout_{\beta_2} \in \mathcal{Y}_{\mathrm CV}(Q_k) \,\, \forall \alpha \in [0,1]$; similarly $Q_1, Q_2 \in \mathcal{Q} \implies \alpha Q_1 + (1-\alpha)Q_2 \in \mathcal{Q} \,\, \forall \alpha \in [0,1]$. Furthermore, we emphasize that each element $\simout \in \mathcal{Y}_{\mathrm CV}(Q_k)$ is a real-valued function having domain $\mathcal{X}$ and identified by some $\beta \in \mathbb{R}^d$; and each element $Q \in \mathcal{Q}$ is a probability measure supported on $\mathcal{X}$ and equivalent to $\native$.
The CV-side update 
selects, under the current sampling law \(Q_k\), a member of
\(\mathcal{Y}_{\mathrm CV}(Q_k)\) having smallest second moment. The LR-side update then
selects a sampling law that reduces the second moment of the newly
chosen estimator while penalizing departure from the native measure
\(P_0\). 

\subsection{Existence and Optimality}

A \emph{fixed point} \((\simout^\star,Q^\star)\) of~\eqref{eq:aem-measure} is a CVLR design that returns the same estimator and measure when passed through the CV-side and LR-side updates in~\eqref{eq:aem-measure}, that is, 
$$\simout^\star
\in
\mathop{\arg\min}_{\simout\in \mathcal{Y}_{\mathrm CV}(Q^\star)} \,\, 
\int \simout(x,Q^\star)^2\,Q^\star(dx),$$
and $$Q^\star
\in
\mathop{\arg\min}_{Q\in\mathcal Q}
\left\{
\int \simout^\star(x,Q)^2\,Q(dx)
+
\lambda D_{\mathrm{KL}}(Q \parallel \native )
\right\}.$$
Thus \(\simout^\star\) is optimal within the chosen estimator class under
\(Q^\star\), while \(Q^\star\) is optimal, subject to regularization
toward \(\native\), for the estimator kernel \(\simout^\star\). The following theorem then identifies conditions under which the design $(\simout^\star,Q^\star)$ is an optimal solution to~\eqref{vrmaster}.

\begin{theorem}[Existence and optimality of AU-CVLR fixed points]
\label{thm:aemm-optimality}
Let
\[
\widetilde g(x):=g(x)-\gamma,
\qquad
R_\beta(x):=\simout_0(x)-\beta^\top\widetilde g(x),
\]
and suppose that
\begin{enumerate}
    \item $\lambda>0$;
    \item $\simout_0,g_1,\ldots,g_d\in L^2(\native)$, i.e., these quantities have finite second moments; and
    \item the centered controls are linearly independent in
    $L^2(\native)$, equivalently,
    \[
    \Sigma_g
    :=
    \int
    \widetilde g(x)\widetilde g(x)^\top
    \,\native(dx)
    \succ0.
    \]
\end{enumerate}
Then the regularized CVLR master problem
\[
\inf_{\beta\in\mathbb R^d,\;Q\in\mathcal Q}
\left\{
\int \simout_\beta(x,Q)^2\,Q(dx)
+
\lambda D_{\rm KL}(Q\|\native)
\right\}
\tag{VR-master}
\label{eq:vrmaster-cvlr}
\]
admits a minimizer $(\beta^\star,Q^\star)$.

Moreover, every fixed point $(\beta^\star,Q^\star)$ of~\eqref{eq:aem-measure}
is a global solution of~\eqref{eq:vrmaster-cvlr}.
Conversely, every global solution of~\eqref{eq:vrmaster-cvlr}
is a fixed point of~\eqref{eq:aem-measure}.  Consequently, the associated CVLR design
kernel
\[
\simout^\star(\cdot,Q)
:=
\simout_{\beta^\star}(\cdot,Q)
\]
and sampling law $Q^\star$ constitute an optimal estimator--measure
design for~\eqref{eq:vrmaster-cvlr}.

Under the stated linear-independence condition, the optimal pair
$(\beta^\star,Q^\star)$ is unique.
\end{theorem}

\begin{proof}
For any $Q\in\mathcal Q$, write
\[
r:=\frac{dQ}{d\native}.
\]
Since $Q$ and $\native$ are equivalent,
$r>0$ $\native$-almost surely and
\[
\int r\,d\native=1.
\]
The CVLR kernel satisfies
\[
\simout_\beta(x,Q)
=
\frac{R_\beta(x)}{r(x)},
\]
and hence the objective in~\eqref{eq:vrmaster-cvlr} can be written as
\[
\Phi_\lambda(\beta,r)
=
\int
\frac{R_\beta(x)^2}{r(x)}
\,\native(dx)
+
\lambda
\int r(x)\log r(x)\,\native(dx),
\tag{22}
\label{eq:phi-beta-r}
\]
over
\[
\beta\in\mathbb R^d,
\qquad
r>0\quad\native\text{-a.s.},
\qquad
\int r\,d\native=1.
\] For each $x$, the map
\[
(z,r)\longmapsto \frac{z^2}{r},
\qquad r>0,
\]
is the perspective of the convex function $z\mapsto z^2$ and is
therefore jointly convex.  Since $R_\beta(x)$ is affine in $\beta$,
the first term in~\eqref{eq:phi-beta-r} is jointly convex in
$(\beta,r)$.  The relative-entropy term is convex in $r$, and is
strictly convex for $\lambda>0$.  Thus $\Phi_\lambda$ is jointly
convex in $\beta$ and $r$.

We first establish existence.  The function $\Phi_{\lambda}$ is proper because it is bounded from below by $-e^{-1}\lambda$, and $\Phi_{\lambda}(0,1) < \infty$ since we have assumed that $\simout_0\in L^2(\native)$. We thus see that~\eqref{vrmaster} has a finite infimum. Now consider a minimizing sequence $(\beta_n,r_n)$, i.e., a sequence $((\beta_n, r_n): n \ge 1)$ such that $\Phi_\lambda(\beta_n, r_n) \to \Phi^*_\lambda$, where $\Phi^*_\lambda$ gives the (finite) infimum.  Since
\[
\lambda
\int r_n\log r_n\,d\native
\]
is bounded above along this sequence, the densities $\{r_n\}$ have
uniformly bounded relative entropy and are therefore uniformly
integrable \cite[\S13.3]{williams1991probability}. Consequently,
after passing to a subsequence,
$r_n$ converges weakly in $L^1(\native)$ to some density $r^\star$
satisfying
\[
r^\star\ge0,
\qquad
\int r^\star\,d\native=1.
\]

The sequence $\{\beta_n\}$ is also bounded.  Indeed, by
Cauchy--Schwarz,
\[
\int
\frac{R_\beta^2}{r}\,d\native
\ge
\left(
\int |R_\beta|\,d\native
\right)^2.
\]
The positive definiteness of $\Sigma_g$ implies the existence of a
constant $c>0$ such that
\[
\int
|u^\top\widetilde g(x)|
\,\native(dx)
\ge c
\qquad
\text{for every }u\in\mathbb R^d,\ \|u\|=1.
\]
Therefore
\[
\int |R_\beta|\,d\native
\ge
c\|\beta\|
-
\int |\simout_0|\,d\native,
\]
so bounded objective values imply bounded $\|\beta\|$.
Passing to a further subsequence gives
$\beta_n\to\beta^\star$.

The perspective functional in~\eqref{eq:phi-beta-r} and the entropy
functional are lower semicontinuous under these convergences.
Therefore $(\beta^\star,r^\star)$ minimizes the relaxed problem with
$r\ge0$.

It remains to note that the minimizer cannot assign zero density to a
set of positive $\native$-measure.  If $r^\star=0$ on such a set,
finiteness of the first term forces
$R_{\beta^\star}=0$ there.  Transferring an arbitrarily small amount
of probability mass to this set changes the first term only to first
order, whereas the entropy contribution contains an
$\epsilon\log\epsilon$ decrease.  For sufficiently small $\epsilon$
this strictly lowers the objective, contradicting optimality.
Hence
\[
r^\star>0
\qquad\native\text{-a.s.},
\]
so the corresponding $Q^\star$ belongs to $\mathcal Q$.  This proves
existence.

We next establish optimality of an AU-CVLR fixed point.  Suppose
$(\beta^\star,Q^\star)$ is a fixed point, and write
\[
r^\star=\frac{dQ^\star}{d\native}.
\]
Optimality of the estimator block gives the normal equations
\[
\int
\widetilde g(x)
\frac{R_{\beta^\star}(x)}
     {r^\star(x)}
\,\native(dx)
=0.
\tag{23}
\label{eq:beta-kkt}
\]
Optimality of the measure block, together with
$\int r\,d\native=1$, yields a scalar multiplier $\nu^\star$ such
that
\[
-
\frac{R_{\beta^\star}(x)^2}
     {r^\star(x)^2}
+
\lambda\bigl(1+\log r^\star(x)\bigr)
+
\nu^\star
=
0,
\qquad
\native\text{-a.s.}
\tag{24}
\label{eq:q-kkt}
\]
Equations~\eqref{eq:beta-kkt}--\eqref{eq:q-kkt}, together with
$\int r^\star\,d\native=1$, are precisely the first-order optimality
conditions for the jointly convex problem~\eqref{eq:phi-beta-r}.
They are therefore sufficient for global optimality.  Hence every
AU-CVLR fixed point solves~\eqref{eq:vrmaster-cvlr}.

Conversely, a global minimizer necessarily minimizes the objective
with respect to each block while the other block is held fixed, and
is therefore an AU-CVLR fixed point.

Finally, strict convexity of the entropy term implies uniqueness of
the optimal density $r^\star$.  Conditional on $r^\star$, the
estimator-side objective is a quadratic in $\beta$ with Hessian
\[
2\int
\frac{\widetilde g(x)\widetilde g(x)^\top}
     {r^\star(x)}
\,\native(dx),
\]
which is positive definite under Assumption~3.  Hence
$\beta^\star$ is also unique.
\end{proof}

\section{A Specific AU-CVLR Algorithm}\label{sec:semaucvlr}

The~\eqref{eq:aem-measure} scheme is conceptual as written because it requires optimization over $\beta \in \mathbb{R}^d$ (CV-side update) and $Q \in \mathcal{Q}$ (LR-side update), which may be difficult. We now propose a specific implementable algorithm to perform these updates. In particular, we show that the CV-side update is ``closed-form'' and reminiscent of the expression for the slope vector in linear regression using least squares, as usual with control variates. On the other hand, the LR-side update is more nuanced but becomes tractable after linearizing the objective appearing in the expression for $Q_{k+1}$. While the LR-side update algorithm we present in Section~\ref{sec:lrside} is nonparametric, further stipulations on the feasible region $\mathcal{Q}$ could simplify the update algorithm, as we demonstrate using an option pricing example in Section~\ref{sec:barrier-asian}. The CV-side update and the LR-side update together provide an implementable nonparametric algorithm algorithm that returns a sequence $(\simout_k, Q_k)$ converging to a fixed point $(\simout^\star,Q^\star)$ of~\eqref{eq:aem-measure}.

\subsection{CV-side Update}\label{sec:cvside}
Denoting $S_k(x):= \frac{dQ_k}{d\native}(x)$ and $L_k(x) := \frac{d\native}{dQ_k}(x) = 1/S_k(x)$, and noticing that each element of $\mathcal{Y}_{\mathrm{CV}}$ is identified by a $\beta \in \mathbb{R}^d$, the CV-side update of~\eqref{eq:aem-measure} becomes 
\begin{align}\label{cvebeta} \beta_{k+1} &:= \mathop{\arg\min}_{\beta \in \mathbb{R}^m} \int R^2_{\beta}(x)L^2_k(x) Q_k(dx) \nonumber \\ &= \mathop{\arg\min}_{\beta \in \mathbb{R}^m} \int R^2_{\beta}(x)L_k(x) \native(dx).
\end{align}
As in the classical linear regression with least-squares loss, the solution $\beta_{k+1}$ to the optimization problem~\eqref{cvebeta} can be obtained through direct differentiation: \begin{equation}\label{betasoln} \beta_{k+1} = \Sigma_k^{-1}\, b_k,\end{equation} where
\begin{align*}
\Sigma_k & := \int (g(x) - \gamma)(g(x) - \gamma)^{\top} \,L_k(x) \native(dx),\\
b_k & :=  \int \simout_0(x) (g(x) - \gamma)\,L_k(x) \, \native(dx).
\end{align*}
(The expression in~\eqref{betasoln} assumes invertibility of $\Sigma_k$, which we implicitly stipulate in the theorem that follows.) The expression in~\eqref{betasoln} plugged into $\simout_{\beta}(\cdot,Q_{k})$ %
then gives the estimator-side update: \begin{equation}\label{estsideupdate} \simout_{k+1} = \frac{\left(\simout_0(x) - b_k^\top \Sigma_k^{-1}(g(x)-\gamma)\right)}{S_k(x)}, \quad S_k(x) = \frac{1}{L_k(x)} \quad x \in \mathcal{X}.
\end{equation}

\subsection{LR-side Update}\label{sec:lrside}

To obtain a measure-side update, we now choose $\mathcal{Q} = \mathcal{P}(\mathcal{X})$ and plug the expression~\eqref{estsideupdate} into~\eqref{eq:aem-measure}: \begin{equation}\label{meassideinitial}
Q_{k+1} \in \mathop{\arg\min}_{Q \in \mathcal{Q}} \left\{ \int R^2_{k+1}(x)\left(\frac{d\native}{dQ}(x)\right)^2 Q(dx) + \lambda D_{\mathrm{KL}}(Q \parallel \native) \right\},
\end{equation}
where we denote the residual $R_{k+1}(x) = \simout_0(x) - \beta_{k+1}^\top(g(x) - \gamma)$. Now assume $\native$ and $Q_k$ have densities $p_0$ and $q$, respectively, with respect to a measure $\nu$, and write
$$\native(dx) = p_0(x)\nu(dx), \quad Q(dx) = q(x)\nu(dx).$$
The optimization problem in~\eqref{meassideinitial} then becomes \begin{equation}\label{meassidedensity}
q_{k+1} \in \mathop{\arg\min}_{q \geq 0 \int q \, = \, 1} \left\{\underbrace{\int R^2_{k+1}(x) \frac{p_0^2(x)}{q(x)} \nu(dx)}_{V_{k+1}(q)} + \lambda \int q(x) \log \frac{q(x)}{p_0(x)} \nu(dx)\right\}.
\end{equation}
Noticing that much of the complexity in the optimization problem in~\eqref{meassidedensity} is due to the (nonlinear in $q$) first term inside the braces, we linearize this term in an attempt at approximate optimization. Toward this, we compute its (von Mises) derivative~\citep{2025yuhenpas} at the density $q_k$ of $Q_{k}$: \begin{align}\label{linapp} DV_{k+1}(q_k)[q-q_k] &= \int -R^2_{k+1}(x) \frac{p_0^2(x)}{q^2_k(x)}(q(x)-q_k(x)) \nu(dx) \nonumber \\ &= \int \underbrace{-R^2_{k+1}(x) \frac{p_0^2(x)}{q^2_k(x)}}_{\psi_{k+1}(x)} Q(dx) + c,\end{align} revealing the influence function of $V_{k+1}$ at $q_k$: $$\psi_{k+1}(x):= -R^2_{k+1}(x)p_0^2(x)/q_k^2(x) = -R^2_{k+1}L_k^2(x).$$ To get a tractable alternative to~\eqref{meassideinitial}, we use~\eqref{linapp} as the linear approximation of the first term in~\eqref{meassideinitial} (after ignoring constants) along with a regularization term that prevents the $(k+1)$-th update $Q_{k+1}$ straying too far from $Q_k$: 
\begin{equation}\label{meassidefinal} Q_{k+1} \in \mathop{\arg\min}_{Q \in \mathcal{P}(\mathcal{X})} \left\{ \int \psi_{k+1}(x) Q(dx) + \lambda D_{\mathrm{KL}}(Q \parallel \native) + \frac{1}{\alpha_k} D_{\mathrm{KL}}(Q \parallel Q_k)\right\},\end{equation} where $\alpha_k>0$ is a step size. Denoting $s(x):= \frac{dQ}{dP_0}(x)$, a further refinement of~\eqref{meassidefinal} gives:
\begin{align}
\min_{\substack{s\geq0\\ \int s\,dP_0=1}}
\left\{
\int \psi_{k+1}(x)s(x)\,P_0(dx)
+
\lambda\int s(x)\log s(x)\,P_0(dx)
\right.\notag \\
+
\left.\frac{1}{\alpha_k}
\int s(x)\left(\log s(x) - \log S_k(x)\right)\,P_0(dx)
\right\}.\label{meassidesoln}
\end{align}
Now explicitly solve~\eqref{meassidesoln} to obtain the LR-side update:
\begin{align}
S_{k+1}(x)
& =
Z_k^{-1}\left[S_k(x) \right]^{\,1/(1+\alpha_k\lambda)}
\exp\!\left\{
\dfrac{\alpha_k}{1+\alpha_k\lambda}
R^2_{k+1}L_k^2(x)
\right\};\notag\\
L_{k+1}(x) &= \frac{1}{S_{k+1}(x)},
\label{eq:aemm-md-closed-form}\end{align}
where the normalizing constant
\[
Z_k
=
\int
\left[S_k(x) \right]^{\,1/(1+\alpha_k\lambda)}
\exp\!\left\{
\frac{\alpha_k}{1+\alpha_k\lambda}
R^2_{k+1}L_k^2(x)
\right\}
P_0(dx).
\]

Expression~\eqref{eq:aemm-md-closed-form} exposes the two forces acting on the LR-side update. The exponential term moves the sampling law in a direction that locally decreases the second moment, while the exponent $(1+\alpha_k\lambda)^{-1}$ shrinks the influence of the current law toward the native law. In particular, when $\lambda=0$, $\frac{dQ_{k+1}}{dQ_k}(x)
\propto
\exp\{-\alpha_k\psi_{k+1}(x)\},$
which is the usual KL mirror-descent update. For $\lambda>0$, the update additionally reflects the regularization toward $P_0$ inherited from the VR-master problem. 

\subsection{Implementation} 
The CV-side and LR-side updates in~\eqref{estsideupdate} and~\eqref{eq:aemm-md-closed-form} suggest a simulation-based algorithm for constructing a variance-reduced estimator of $\theta$. Specifically, a one-sample estimator $\theta$ can be constructed through the following three steps:
\begin{tcolorbox}[
    enhanced,
    colback=white,
    colframe=black,
    boxrule=0.6pt,
    arc=0pt,
    outer arc=0pt,
    attach boxed title to top left={
        xshift=2mm,
        yshift=-1.5mm
    },
    boxed title style={
        colback=white,
        colframe=white,
        boxrule=0pt
    },
    title={Nonparametric AU-CVLR},
    fonttitle=\bfseries,
    coltitle=black,
    left=2mm,
    right=2mm,
    top=2mm,
    bottom=1mm
]
\begin{enumerate} \item[(a)] starting with an initial estimate $\hat{S}_0(x) = 1$ of the likelihood ratio $\frac{dQ_k}{d\native}(x)$, and the ability to generate from the native law $\native$, execute the updates in~\eqref{estsideupdate} and~\eqref{eq:aemm-md-closed-form} for (say) $k$ iterations to obtain the estimator $\simout_k$ and the likelihood ratio $S_k$ (see Algorithm~\ref{alg:aemm-md-mcmc}); \item[(b)] use \emph{rejection sampling}~\citep[\S6.4.1]{nelson2013foundations} (see Algorithm~\ref{alg:rejection-qk}) or MCMC with Metropolis-Hastings~\citep{chib1995understanding} to generate $X \sim Q_k$ (see Algorithm~\ref{alg:mcmc-qk}); \item[(c)] Construct $\simout_{k}(X,Q_k), X \sim Q_k$.
\end{enumerate}
\end{tcolorbox}

\begin{algorithm}[htb]
\caption{Alternating Updates CVLR Design Kernel}
\label{alg:aemm-md-mcmc}
\begin{algorithmic}[1]

\Require Native sampler $P_0$; simulation oracle $\simout_0(x)$;
controls $g$ with known mean $\gamma=\E_{P_0}[g(X)]$;
regularization parameter $\lambda\geq 0$;
stepsizes $\{\alpha_k\}$;
sample sizes $n_k$ and $m_k$.

\State Define $\widetilde g(x)\gets g(x)-\gamma$.
\State Initialize $Q_0=P_0$, $\widehat S_0(x)\equiv 1$, where
$\widehat S_k(x)\approx \frac{dQ_k}{dP_0}(x).$
\State Generate an initial sample
$X_{0,1},\ldots,X_{0,n_0}\sim P_0,$
and evaluate $\simout_0(X_{0,j})$.

\For{$k=0,1,2,\ldots, K-1$}

    \Statex \textbf{Estimator update}

    \State Set $\widehat L_k(x)\gets \frac{1}{\widehat S_k(x)}
    \approx \frac{dP_0}{dQ_k}(x).$

    \State Using the current sample
    $X_{k,1},\ldots,X_{k,n_k}\approx Q_k$, compute
    \begin{align*}
    \widehat\Sigma_k
    & \gets
    \frac{1}{n_k}
    \sum_{j=1}^{n_k}
    \widetilde g(X_{k,j})
    \widetilde g(X_{k,j})^\top
    \widehat L_k(X_{k,j})^2\\ \widehat b_k
    & \gets
    \frac{1}{n_k}
    \sum_{j=1}^{n_k}
    \widetilde g(X_{k,j})
    \simout_0(X_{k,j})
    \widehat L_k(X_{k,j})^2 .
    \end{align*}

    \State Compute $\widehat\beta_{k+1}
    \gets
    \widehat\Sigma_k^{-1}\widehat b_k.$

    \State Define the estimated measure-side influence function
    \[
    \widehat\psi_{k+1}(x)
    \gets
    -
    \left[
    \left\{
    \simout_0(x)
    -
    \widehat\beta_{k+1}^{\top}\widetilde g(x)
    \right\}
    \widehat L_k(x)
    \right]^2.
    \]

    \Statex \textbf{Measure update}

    \State Set $a_k\gets\frac{1}{1+\alpha_k\lambda},
    \qquad
    c_k\gets\frac{\alpha_k}{1+\alpha_k\lambda}.$

    \State Define the unnormalized Radon--Nikodym derivative
    \[
    \overline S_{k+1}(x)
    \gets
    \widehat S_k(x)^{a_k}
    \exp\!\left\{-c_k\widehat\psi_{k+1}(x)\right\}.
    \]

    \State Generate $U_{k,1},\ldots,U_{k,m_k}\stackrel{\mathrm{iid}}{\sim}P_0$
    and evaluate $\simout_0(U_{k,\ell})$ to compute
    $\overline S_{k+1}(U_{k,\ell})$.

    \State Estimate the normalizing constant $\widehat Z_k
    \gets
    \frac{1}{m_k}
    \sum_{\ell=1}^{m_k}
    \overline S_{k+1}(U_{k,\ell}).$

    \State Set $\widehat S_{k+1}(x)
    \gets
    \frac{\overline S_{k+1}(x)}
         {\widehat Z_k}.$

\EndFor

\end{algorithmic}
\end{algorithm}

It is important that the updates in~\eqref{estsideupdate} and~\eqref{eq:aemm-md-closed-form} can be carried out with only the likelihood ratios and the ability to sample from the native law $\native$. The property that the sampling law $Q_k$ participates in the recursive updates only through the likelihood ratios is inherited from the product form of $\simout_{\beta}$,
and seems to contribute crucially to stability during implementation. If rejection sampling is used in (b), the accepted observation is exactly distributed according to $Q_k$ without having to know the normalizing constant that relates $\overline S_k$ to $dQ_k/dP_0$.

\begin{algorithm}[htb]
\caption{Rejection Sampling from $Q_k$}
\label{alg:rejection-qk}
\begin{algorithmic}[1]

\Require Native sampler $P_0$; unnormalized density ratio
$\overline S_k:\X\to(0,\infty)$ satisfying
\[
Q_k(dx)\propto \overline S_k(x)\,P_0(dx);
\]
constant $M_k<\infty$ such that
\[
\overline S_k(x)\leq M_k
\qquad\text{for all }x\in\X.
\]

\Repeat
    \State Generate $X'\sim P_0$.
    \State Generate $U\sim\mathrm{Unif}(0,1)$ independently.
\Until{$U\leq \overline S_k(X')/M_k$}

\State \Return $X\gets X'$.

\end{algorithmic}
\end{algorithm}

\begin{algorithm}[t]
\caption{Independence Metropolis--Hastings Sampling from $Q_k$}
\label{alg:mcmc-qk}
\begin{algorithmic}[1]

\Require Native sampler $P_0$; unnormalized density ratio
$\overline r_k:\X\to(0,\infty)$ satisfying
\[
Q_k(dx)\propto \overline r_k(x)\,P_0(dx);
\]
initial state $X^{(0)}$; burn-in length $B$; desired sample size $N$.

\For{$\ell=1,\ldots,B+N$}

    \State Generate a proposal
    \[
    X'\sim P_0.
    \]

    \State Compute the acceptance probability
    \[
    a\bigl(X^{(\ell-1)},X'\bigr)
    \gets
    1\wedge
    \frac{\overline r_k(X')}
         {\overline r_k(X^{(\ell-1)})}.
    \]

    \State Generate $U\sim\mathrm{Unif}(0,1)$.

    \If{$U\leq a\bigl(X^{(\ell-1)},X'\bigr)$}
        \State $X^{(\ell)}\gets X'$.
    \Else
        \State $X^{(\ell)}\gets X^{(\ell-1)}$.
    \EndIf

\EndFor

\State \Return
\[
X^{(B+1)},\ldots,X^{(B+N)}
\]
as an MCMC sample from $Q_k$.

\end{algorithmic}
\end{algorithm}

\begin{theorem}[Convergence of Nonparametric AU-CVLR]
\label{thm:aemm-md-convergence}
Let
\[
\widetilde g(x):=g(x)-\gamma,
\qquad
R_\beta(x):=\simout_0(x)-\beta^\top\widetilde g(x),
\]
and, for a density $S=dQ/d\native$, define the 
second-moment functional
\[
F(S)
:=
\min_{\beta\in\mathbb R^d}
\int
\frac{R_\beta(x)^2}{S(x)}
\,\native(dx).
\]
Let
\[
H(S):=\int S(x)\log S(x)\,\native(dx)
\]
and
\[
\Psi_\lambda(S):=F(S)+\lambda H(S),
\qquad \lambda>0.
\]

Suppose the following conditions hold.

\begin{enumerate}
\item[(A1)]
The regularized CVLR problem has a unique minimizer
$(\beta^\star,S^\star)$, with
\[
S^\star>0\quad \native\text{-a.s.},
\qquad
\int S^\star\,d\native=1.
\]

\item[(A2)]
For every density $S$ visited by the algorithm, the estimator-side
problem defining $F(S)$ has a unique minimizer, denoted $\beta(S)$,
and the map
\[
S\longmapsto\beta(S)
\]
is continuous in $L^1(\native)$.

\item[(A3)]
The functional $F$ is differentiable on a convex set
$\mathcal R$ containing $S^\star$ and all iterates, and is
$L$-smooth relative to $H$ on $\mathcal R$, i.e.,
\[
F(S')
\le
F(S)
+
\int
\psi_S(x)\{S'(x)-S(x)\}\,\native(dx)
+
L D_{\rm KL}(Q'\|Q)
\]
for every $S,S'\in\mathcal R$, where
\[
\psi_S(x)
=
-\frac{R_{\beta(S)}(x)^2}{S(x)^2}
\]
is the first-variation density of $F$ at $S$ and
$Q,Q'$ are the probability measures having densities $S,S'$
with respect to $\native$.

\item[(A4)]
The mirror-descent stepsizes satisfy
\[
0<\underline\alpha
\le
\alpha_k
\le
\frac{1}{L},
\qquad k\ge0.
\]
\end{enumerate}

Starting from $S_0\in\mathcal R$, let
\[
\beta_{k+1}:=\beta(S_k)
\]
be the exact estimator-side update, and let $S_{k+1}$ be given by
the AU-CVLR mirror descent update
\begin{align}
S_{k+1}
& \in
\arg\min_{\substack{S\in\mathcal R\\ \int S\,d\native=1}}
\left\{
\int
\psi_{k+1}(x)\{S(x)-S_k(x)\}\,\native(dx)
+
\lambda H(S)\right. \notag\\
& \quad +
\left.\frac{1}{\alpha_k}
D_{\rm KL}(Q\|Q_k)
\right\},
\label{eq:aemm-md-prox}
\end{align}
where
\[
\psi_{k+1}(x)
=
-\frac{R_{\beta_{k+1}}(x)^2}{S_k(x)^2}.
\]
Equivalently, $S_{k+1}$ is given by the closed-form update
in~\eqref{eq:aemm-md-closed-form}.

Then
\[
\Psi_\lambda(S_k)\longrightarrow
\Psi_\lambda(S^\star),
\]
and
\[
D_{\rm KL}(Q_k\|Q^\star)\longrightarrow0.
\]
Consequently,
\[
\|S_k-S^\star\|_{L^1(\native)}
\longrightarrow0,
\]
and
\[
\beta_{k+1}\longrightarrow\beta^\star.
\]
Hence
\[
(Q_k,\simout_{\beta_{k+1}}(\cdot,Q_k))
\longrightarrow
(Q^\star,\simout_{\beta^\star}(\cdot,Q^\star)),
\]
and the limiting pair is the unique fixed point of AU-CVLR mirror descent and,
by Theorem~\ref{thm:aemm-optimality}, the unique optimal CVLR
design for~\eqref{vrmaster}.
\end{theorem}

\begin{proof}
The measure update~\eqref{eq:aemm-md-prox} is a composite
KL mirror-descent step for
\[
\Psi_\lambda=F+\lambda H.
\]
The three-point inequality for Bregman proximal updates gives, for
every feasible density $S$,
\begin{align}
&\alpha_k
\int
\psi_{k+1}(x)
\{S_{k+1}(x)-S(x)\}\,\native(dx)
\nonumber\\
&\qquad
+
\alpha_k\lambda
\{H(S_{k+1})-H(S)\}
\nonumber\\
&\le
D_{\rm KL}(Q\|Q_k)
-
D_{\rm KL}(Q\|Q_{k+1})
-
D_{\rm KL}(Q_{k+1}\|Q_k).
\tag{26}
\label{eq:three-point}
\end{align}

Set $S=S^\star$.  Relative smoothness of $F$ gives
\[
F(S_{k+1})
\le
F(S_k)
+
\langle\psi_{k+1},S_{k+1}-S_k\rangle
+
L D_{\rm KL}(Q_{k+1}\|Q_k),
\]
whereas convexity of $F$ gives
\[
F(S^\star)
\ge
F(S_k)
+
\langle\psi_{k+1},S^\star-S_k\rangle.
\]
Subtracting these inequalities yields
\[
F(S_{k+1})-F(S^\star)
\le
\langle\psi_{k+1},S_{k+1}-S^\star\rangle
+
L D_{\rm KL}(Q_{k+1}\|Q_k).
\]
Similarly, convexity of $H$ implies
\[
H(S_{k+1})-H(S^\star)
\le
\langle\nabla H(S_{k+1}),
S_{k+1}-S^\star\rangle.
\]
Combining these relations with~\eqref{eq:three-point} gives
\begin{align}
\alpha_k
\{\Psi_\lambda(S_{k+1})-\Psi_\lambda(S^\star)\}
&\le
D_{\rm KL}(Q^\star\|Q_k)
-
D_{\rm KL}(Q^\star\|Q_{k+1})
\nonumber\\
&\quad
-
(1-\alpha_kL)
D_{\rm KL}(Q_{k+1}\|Q_k).
\tag{27}
\label{eq:aemm-descent}
\end{align}
Since $\alpha_kL\le1$, the final term is nonpositive and may be
discarded.  Summing~\eqref{eq:aemm-descent} from $k=0$ to $N-1$
therefore yields
\[
\sum_{k=0}^{N-1}
\alpha_k
\{\Psi_\lambda(S_{k+1})-\Psi_\lambda(S^\star)\}
\le
D_{\rm KL}(Q^\star\|Q_0).
\]
Because $\alpha_k\ge\underline\alpha>0$ and every summand is
nonnegative,
\[
\Psi_\lambda(S_k)-\Psi_\lambda(S^\star)
\longrightarrow0.
\]

Moreover, $\Psi_\lambda$ is $\lambda$-strongly convex relative to
$H$.  Indeed, $F$ is convex and the term $\lambda H$ contributes
$\lambda D_{\rm KL}(\cdot\|\cdot)$.  At the minimizer $S^\star$,
the first-order optimality condition and the normalization constraint
therefore imply
\[
\Psi_\lambda(S)-\Psi_\lambda(S^\star)
\ge
\lambda
D_{\rm KL}(Q\|Q^\star)
\]
for every feasible $S$.  Hence
\[
D_{\rm KL}(Q_k\|Q^\star)\longrightarrow0.
\]
Pinsker's inequality then gives
\[
\|Q_k-Q^\star\|_{\rm TV}\longrightarrow0,
\]
or equivalently
\[
\|S_k-S^\star\|_{L^1(\native)}\longrightarrow0.
\]

Finally, Assumption~(A2) yields
\[
\beta_{k+1}
=
\beta(S_k)
\longrightarrow
\beta(S^\star)
=
\beta^\star.
\]
Thus the estimator and measure iterates converge to
$(\beta^\star,Q^\star)$, which is the unique fixed point and the
unique global CVLR optimum.
\end{proof}

\section{NUMERICAL ILLUSTRATION}
\label{sec:numerical}
We illustrate CVLR estimators and their identification using three examples. In each case, we compare the proposed CVLR estimator against CMC, pure-CV and pure-LR estimators in terms of the extent of variance reduction.

\subsection{Finite-horizon M/M/1 queue}
\label{sec:mm1-example}

Consider a finite-horizon \(M/M/1\) queue initialized empty.
Interarrival times \(A_i\) are independent exponential random
variables with rate \(0.9\), service times \(S_i\) are independent
exponential random variables with rate \(1\), and the waiting times
evolve according to the Lindley recursion
\[
W_1=0,
\qquad
W_{i+1}
=
\bigl(W_i+S_i-A_i\bigr)^+,
\qquad
i=1,\ldots,19.
\]
We use throughout the same two control variates,
\[
g_1
=
\frac{1}{20}\sum_{i=1}^{20}S_i,
\qquad
g_2
=
\frac{1}{19}\sum_{i=1}^{19}A_i,
\]
whose native expectations are known exactly:
\[
\E_{\native}[g_1]=1,
\qquad
\E_{\native}[g_2]=\frac{1}{0.9}.
\]
These controls are deliberately low dimensional.  They summarize
the average service and arrival environment experienced by the
simulated trajectory while leaving substantial dependence on the
ordering and interaction of the individual arrival and service
times.

We focus on the exceedance probability
\[
\theta_q
=
\mathbb{P}_{\native}(W_{20}>c_q),
\]
where the numerical thresholds \(c_q\) are calibrated to lie
approximately at selected \(q\)-quantiles of the native waiting-time
distribution.  The thresholds used are
\[
\begin{array}{c|cccccccc}
q
&0.50&0.60&0.70&0.80&0.90&0.95&0.99&0.999\\ \hline
c_q
&2.32&3.25&4.37&5.82&8.06&10.07&14.16&19.19.
\end{array}
\]

\begin{remark}
The thresholds \(c_q\) were calibrated using an independent benchmark
run of \(2\times10^6\) native trajectories.  This run is used only to
define the eight test instances.  Each competing method is
subsequently supplied the numerical threshold \(c_q\), but neither
the corresponding quantile label nor the associated exceedance
probability is used in fitting, tuning, or budget allocation.
Thus \(q\) should be viewed as a convenient label for the severity
of the test instance rather than as information supplied to the
estimator.
\end{remark}

The expanded range of thresholds is useful because it allows the
same stochastic model and the same two controls to exhibit several
different variance-reduction regimes.  At the more central
thresholds, average arrival and service conditions contain
substantial information about whether \(W_{20}\) exceeds \(c_q\),
and control variates are therefore expected to be particularly
effective.  As the threshold moves farther into the tail, the event
is increasingly associated with an unusual region of trajectory
space, so changing the sampling measure becomes progressively more
important.

We compare CMC, pure-CV, pure-LR, and AU-CVLR under a common
simulation budget.  We consider the budgets $n=10{,}000$ and $n=100{,}000$.
For each variance-reduction method, \(10\%\) of the total budget is
reserved as a pilot, that is, the pilot has length $B=0.1n$.
CMC requires no training and uses all \(n\) observations directly,
\[
\widehat\theta_{\rm CMC}
=
\frac{1}{n}
\sum_{i=1}^{n}
\mathbf 1\{W_{20}^{(i)}>c_q\},
\qquad
W_{20}^{(i)}
\stackrel{\mathrm{iid}}{\sim}
\native.
\]
For pure CV, all \(B\) pilot trajectories are used to estimate the
control coefficient \(\beta\).  The remaining $N=n-B=0.9n$
fresh native trajectories are used in the final CV estimator.
Thus observations used to fit the control coefficient are charged
to the method and are not reused in the final sample average.

For pure LR and AU-CVLR, the same \(B\)-observation native pilot is
used to construct and tune the sampling design.  The CV-side and
LR-side quantities appearing in the alternating updates can be
represented as native-measure expectations.  Consequently, an
observation generated from \(\native\) is charged only once even
though its output may be reused over several alternating updates.
This seems appropriate in settings where obtaining simulation
replications is expensive relative to carrying out the optimization
steps.

The pilot is divided into a \(70\%\) fitting subset and a \(30\%\)
validation subset.  Candidate sampling designs are constructed
using the fitting subset.  The validation subset is then used to
select the regularization parameter \(\lambda\) and to estimate the
quantities required to divide the remaining budget between
likelihood-ratio normalization and production.  Once \(\lambda\)
has been selected, the corresponding design is refitted using the
complete \(B\)-observation pilot.  No production observation is
used for fitting, tuning, or model selection.

For every threshold and both budget levels, we use
\[
K=64,
\qquad
\alpha_k\equiv0.2
\]
alternating updates and fixed mirror stepsize, respectively.
Neither \(K\) nor \(\alpha_k\) is tuned separately across thresholds
or budgets.  The candidate grid used to select the regularization
parameter is also fixed in advance:
\[
\begin{split}
\Lambda=\{&
10^{-4},
3\times10^{-4},
10^{-3},
3\times10^{-3},
10^{-2},
3\times10^{-2},\\
&0.1,\,
0.3,\,
1,\,
3,\,
10,\,
30,\,
100
\}.
\end{split}
\]
For AU-CVLR, the pure-CV endpoint
\(\lambda=\infty\), corresponding to no measure adaptation, is
included as an additional candidate.  Thus the algorithm is free
to revert to ordinary CV when the pilot provides insufficient
evidence that a change of measure is useful.

After a finite-\(\lambda\) sampling law has been selected, write its
unnormalized density ratio as \(s\), so that
\[
\frac{dQ}{d\native}(x)
=
\frac{s(x)}{C},
\qquad
C
=
\E_{\native}[s(X)].
\]
The \(0.9n\) observations remaining after the pilot are divided into
\(m\) fresh native observations used to estimate \(C\) and \(N\)
fresh production observations generated according to the learned
sampling law \(Q\), so that $B+m+N=n$.
The resulting AU-CVLR estimator is
\[
\widehat\theta_{\rm CVLR}
=
\frac{1}{N}
\sum_{i=1}^{N}
\left[
\mathbf 1\{W_{20}^{(i)}>c_q\}
-
\widehat\beta^\top
\bigl(g(X_i)-\gamma\bigr)
\right]
\frac{\widehat C}{s(X_i)},
\qquad
X_i\sim\widehat Q.
\]
Pure LR uses the same construction with
\(\widehat\beta=0\).

Although the fitted control coefficient, the selected regularization
parameter, and the learned sampling measure are all random, their
dependence on the pilot does not introduce bias because the
production sample is independent of the pilot.  Conditional on any
fixed pilot realization, the CV correction has native mean zero and
the exact likelihood-ratio estimator has mean \(\theta_q\).
Likewise, replacing the normalizing constant \(C\) by
\[
\widehat C
=
\frac1m\sum_{j=1}^m s(Z_j),
\qquad
Z_j\stackrel{\mathrm{iid}}{\sim}\native,
\]
preserves unconditional unbiasedness because \(\widehat C\) is an
unbiased estimate of \(C\) and is obtained independently of the
production sample.  Consequently, $\E[\widehat\theta]=\theta_q$
for each of the procedures considered here, and hence
$\operatorname{MSE}(\widehat\theta)=
\var(\widehat\theta)$. We therefore report variance throughout this example.%

\begin{remark}
Our simulation-budget accounting treats each retained production
observation from the learned measure \(\widehat Q\) as one simulation
observation.  The learned density ratio
\(\widehat S=d\widehat Q/d\native\) determines \(\widehat Q\), but
draws from this measure may in practice be obtained in several ways,
including rejection sampling or Markov chain Monte Carlo.  We do
not charge additional proposals, rejected trajectories, or Markov
chain transitions associated with a particular implementation of
the sampling step.  Thus the comparisons below measure statistical
efficiency under a common number of usable simulation outputs,
rather than the computational efficiency of a particular mechanism
for sampling from \(\widehat Q\).
\end{remark}

The split between \(m\) and \(N\) is determined exclusively from
pilot information.  If
\[
c_s^2
=
\frac{
\var_{\native}(s(X))
}{
\E_{\native}[s(X)]^2
}
\]
and \(V\) denotes the one-sample variance under the learned sampling
law when the exact normalizing constant is used, then, conditional
on a fixed design, the variance of the estimator using the
independent estimate \(\widehat C\) is
\[
\var(\widehat\theta)
=
\frac{\theta_q^2c_s^2}{m}
+
\frac{V}{N}
\left(
1+\frac{c_s^2}{m}
\right).
\]
The quantities in this expression are replaced by estimates obtained
from the pilot.  In particular, the unknown \(\theta_q\) is
replaced by the pilot exceedance frequency; the benchmark value of
\(\theta_q\) is never used in allocation or tuning.

For a method \(\mathcal M\), let
\[
\widehat V_{n,\mathcal M}
\]
denote the estimated variance of its final estimator under total
simulation budget \(n\), including the randomness associated with
pilot fitting, design selection, normalization, and production.
In parallel with examples in subsequent sections, we report the
variance-reduction factor
\[
\widehat{\mathrm{VRF}}_{n,\mathcal M}
=
\frac{
\widehat V_{n,\rm CMC}
}{
\widehat V_{n,\mathcal M}
}.
\]
Thus \(\widehat{\mathrm{VRF}}_n>1\) favors the corresponding
variance-reduction method over CMC, while a value below one indicates
that, at the specified total simulation budget, the cost of learning
and implementing the variance-reduction design outweighs the
variance reduction it provides.

Whenever a finite-\(\lambda\) design is selected, we impose the
pre-specified bounds
\[
0.005n
\le m
\le
\frac{1}{2}(n-B).
\]
The lower bound prevents an unrealistically small sample from being
used to normalize the likelihood ratio.  The upper bound guarantees
that at least one half of the post-pilot budget remains available
for production and prevents a noisy validation sample from assigning
nearly the entire remaining budget to normalization.

A further safeguard is used for very rare events.  If fewer than
five exceedances are observed in the complete \(B\)-observation
pilot, the pilot is regarded as providing insufficient information
for reliable measure adaptation.  Pure LR then reverts to CMC for
the remaining budget, while AU-CVLR reverts to pure CV.  This rule
is fixed in advance and prevents the procedure from attempting to
learn the location of a rare-event region that is essentially absent
from the training sample.

Table~\ref{tab:mm1-var-10k} reports the results for
\(n=10{,}000\). Panel (a) gives the variance of the final
finite-budget estimator, while Panel (b) gives the corresponding
variance-reduction factor.

\begin{table}[t]
\centering
\small
\caption{Finite-budget variance comparison for
\(n=10{,}000\).  Panel (a) reports variance \(\times10^6\);
Panel (b) reports the variance-reduction factor
\(\widehat{\mathrm{VRF}}_n
=\widehat V_{n,\rm CMC}/\widehat V_n\).
Boldface indicates the best-performing method in each row.}
\label{tab:mm1-var-10k}

\begin{tabular}{@{}crrrr@{}}
\toprule
\multicolumn{5}{c}{(a) Variance \(\times10^6\)}\\
\cmidrule(lr){1-5}
\(q\)
& CMC
& Pure-CV
& Pure-LR
& AU-CVLR\\
\midrule
0.50  & 25.000 & \textbf{16.98} & 27.79 & 17.32\\
0.60  & 24.000 & 15.56 & 25.89 & \textbf{15.50}\\
0.70  & 21.000 & 13.55 & 20.81 & \textbf{13.19}\\
0.80  & 16.003 & 10.908 & 12.893 & \textbf{9.792}\\
0.90  & 8.995 & 7.034 & 5.084 & \textbf{4.634}\\
0.95  & 4.736 & 4.171 & 1.837 & \textbf{1.770}\\
0.99  & 0.980 & 1.012 & 0.299 & \textbf{0.229}\\
0.999 & \textbf{0.0962} & 0.1068 & 0.1068 & 0.1068\\
\bottomrule
\end{tabular}

\vspace{0.7em}

\begin{tabular}{@{}crrrr@{}}
\toprule
\multicolumn{5}{c}{(b) Variance-reduction factor}\\
\cmidrule(lr){1-5}
\(q\)
& CMC
& Pure-CV
& Pure-LR
& AU-CVLR\\
\midrule
0.50  & 1.00 & \textbf{1.47} & 0.90 & 1.44\\
0.60  & 1.00 & 1.54 & 0.93 & \textbf{1.55}\\
0.70  & 1.00 & 1.55 & 1.01 & \textbf{1.59}\\
0.80  & 1.00 & 1.47 & 1.24 & \textbf{1.63}\\
0.90  & 1.00 & 1.28 & 1.77 & \textbf{1.94}\\
0.95  & 1.00 & 1.14 & 2.58 & \textbf{2.68}\\
0.99  & 1.00 & 0.97 & 3.28 & \textbf{4.28}\\
0.999 & \textbf{1.00} & 0.90 & 0.90 & 0.90\\
\bottomrule
\end{tabular}
\end{table}

Representative pilot-selected LR and AU-CVLR tuning parameters for
this budget are shown in Table~\ref{tab:mm1-tuning-10k}.  The
selected values vary across macroreplications because both
\(\lambda\) and the allocation are determined from the pilot.  The
entries summarize typical selections and make explicit the magnitude
of the resulting normalization and production budgets.  In every
row,
\[
m+N=9000,
\]
following the charged \(1000\)-trajectory pilot.

\begin{table}[t]
\centering
\small
\caption{Representative pilot-selected tuning and post-pilot
allocations for \(n=10{,}000\).  The notation \(\infty\) denotes
the pure-CV endpoint.}
\label{tab:mm1-tuning-10k}

\begin{tabular}{@{}crrr@{}}
\toprule
\multicolumn{4}{c}{Pure LR}\\
\cmidrule(lr){1-4}
\(q\)
& \(\lambda\)
& \(m\)
& \(N\)\\
\midrule
0.50  & \(100\)       & 50   & 8950\\
0.60  & \(3\)         & 960  & 8040\\
0.70  & \(1\)         & 1630 & 7370\\
0.80  & \(0.3\)       & 2348 & 6652\\
0.90  & \(0.03\)      & 3268 & 5732\\
0.95  & \(10^{-4}\)   & 2977 & 6023\\
0.99  & \(10^{-4}\)   & 718  & 8282\\
0.999 & fallback      & 0    & 9000\\
\bottomrule
\end{tabular}

\vspace{0.7em}

\begin{tabular}{@{}crrr@{}}
\toprule
\multicolumn{4}{c}{AU-CVLR}\\
\cmidrule(lr){1-4}
\(q\)
& \(\lambda\)
& \(m\)
& \(N\)\\
\midrule
0.50  & \(\infty\)    & 0    & 9000\\
0.60  & \(3\)         & 950  & 8050\\
0.70  & \(1\)         & 1410 & 7590\\
0.80  & \(0.3\)       & 1704 & 7296\\
0.90  & \(0.03\)      & 2933 & 6067\\
0.95  & \(10^{-4}\)   & 2692 & 6308\\
0.99  & \(10^{-4}\)   & 740  & 8260\\
0.999 & fallback      & 0    & 9000\\
\bottomrule
\end{tabular}
\end{table}

The additional central thresholds make the transition between the
two variance-reduction mechanisms more visible.  At the
\(0.50\)-quantile, pure CV is clearly preferable to pure LR:
the controls explain a substantial portion of the event indicator,
whereas changing the measure provides too little additional benefit
to justify its training and normalization costs.  Pure CV achieves
a variance-reduction factor of \(1.47\), while pure LR has a factor
below one.  AU-CVLR correspondingly selects, or remains very close
to, the pure-CV endpoint.  The same qualitative behavior persists
at the \(0.60\)-quantile.

By the \(0.70\)- and \(0.80\)-quantiles, moderate measure adaptation
begins to provide an incremental benefit beyond CV.  At the
\(0.90\)-quantile, the ordering changes substantially: the
variance-reduction factors are \(1.28\), \(1.77\), and \(1.94\) for
pure CV, pure LR, and AU-CVLR, respectively.  Thus pure LR now
improves clearly on pure CV, while their joint use improves further
on both.  At the \(0.95\)- and \(0.99\)-quantiles, the balance shifts
still farther toward measure adaptation and the incremental
contribution of the controls becomes small.  At \(q=0.99\),
AU-CVLR achieves a variance-reduction factor of approximately
\(4.28\).

The \(0.999\)-quantile behaves differently.  With \(B=1000\), the
pilot contains only about one exceedance on average.  The
five-exceedance safeguard therefore triggers in essentially all
replications.  Pure LR and AU-CVLR spend \(1000\) observations
determining that reliable measure adaptation is unsupported by the
pilot and then use only the remaining \(9000\) observations for
estimation.  CMC, in contrast, uses all \(10{,}000\) observations
directly.  Its variance-reduction factor is therefore one, compared
with approximately \(0.90\) for each of the variance-reduction
procedures.  This is a consequence of finite training information
rather than evidence that native sampling is intrinsically
preferable for the extreme-tail problem.

Table~\ref{tab:mm1-var-100k} reports the corresponding results for
\(n=100{,}000\).  Panel (a) gives variance multiplied by \(10^7\),
and Panel (b) gives the corresponding variance-reduction factors.

\begin{table}[t]
\centering
\small
\caption{Finite-budget variance comparison for
\(n=100{,}000\).  Panel (a) reports variance \(\times10^7\);
Panel (b) reports the corresponding variance-reduction factor.
Boldface indicates the best-performing method in each row.}
\label{tab:mm1-var-100k}

\begin{tabular}{@{}crrrr@{}}
\toprule
\multicolumn{5}{c}{(a) Variance \(\times10^7\)}\\
\cmidrule(lr){1-5}
\(q\)
& CMC
& Pure-CV
& Pure-LR
& AU-CVLR\\
\midrule
0.50  & 25.000 & \textbf{16.96} & 27.78 & 16.98\\
0.60  & 24.000 & \textbf{15.54} & 25.82 & 15.59\\
0.70  & 21.000 & 13.52 & 20.64 & \textbf{13.06}\\
0.80  & 16.003 & 10.866 & 13.008 & \textbf{9.777}\\
0.90  & 8.995 & 7.015 & 5.177 & \textbf{4.660}\\
0.95  & 4.736 & 4.161 & 1.880 & \textbf{1.817}\\
0.99  & 0.980 & 1.007 & \textbf{0.236} & 0.237\\
0.999 & 0.0962 & 0.1055 & 0.0305 & \textbf{0.0218}\\
\bottomrule
\end{tabular}

\vspace{0.7em}

\begin{tabular}{@{}crrrr@{}}
\toprule
\multicolumn{5}{c}{(b) Variance-reduction factor}\\
\cmidrule(lr){1-5}
\(q\)
& CMC
& Pure-CV
& Pure-LR
& AU-CVLR\\
\midrule
0.50  & 1.00 & \textbf{1.47} & 0.90 & 1.47\\
0.60  & 1.00 & \textbf{1.54} & 0.93 & 1.54\\
0.70  & 1.00 & 1.55 & 1.02 & \textbf{1.61}\\
0.80  & 1.00 & 1.47 & 1.23 & \textbf{1.64}\\
0.90  & 1.00 & 1.28 & 1.74 & \textbf{1.93}\\
0.95  & 1.00 & 1.14 & 2.52 & \textbf{2.61}\\
0.99  & 1.00 & 0.97 & \textbf{4.15} & 4.14\\
0.999 & 1.00 & 0.91 & 3.15 & \textbf{4.41}\\
\bottomrule
\end{tabular}
\end{table}

The corresponding representative pilot-selected tuning parameters
and post-pilot allocations are shown in
Table~\ref{tab:mm1-tuning-100k}.  Here
\[
m+N=90{,}000,
\]
following the charged \(10{,}000\)-trajectory pilot.

\begin{table}[t]
\centering
\small
\caption{Representative pilot-selected tuning and post-pilot
allocations for \(n=100{,}000\).}
\label{tab:mm1-tuning-100k}

\begin{tabular}{@{}crrr@{}}
\toprule
\multicolumn{4}{c}{Pure LR}\\
\cmidrule(lr){1-4}
\(q\)
& \(\lambda\)
& \(m\)
& \(N\)\\
\midrule
0.50  & \(100\)       & 500   & 89500\\
0.60  & \(3\)         & 9330  & 80670\\
0.70  & \(1\)         & 16100 & 73900\\
0.80  & \(0.3\)       & 24013 & 65987\\
0.90  & \(0.03\)      & 31927 & 58073\\
0.95  & \(10^{-4}\)   & 29421 & 60579\\
0.99  & \(10^{-4}\)   & 7578  & 82422\\
0.999 & \(10^{-4}\)   & 818   & 89182\\
\bottomrule
\end{tabular}

\vspace{0.7em}

\begin{tabular}{@{}crrr@{}}
\toprule
\multicolumn{4}{c}{AU-CVLR}\\
\cmidrule(lr){1-4}
\(q\)
& \(\lambda\)
& \(m\)
& \(N\)\\
\midrule
0.50  & \(10\)        & 4560  & 85440\\
0.60  & \(3\)         & 8960  & 81040\\
0.70  & \(1\)         & 13600 & 76400\\
0.80  & \(0.3\)       & 19897 & 70103\\
0.90  & \(0.1\)       & 22942 & 67058\\
0.95  & \(10^{-4}\)   & 26703 & 63297\\
0.99  & \(10^{-4}\)   & 7294  & 82706\\
0.999 & \(10^{-4}\)   & 877   & 89123\\
\bottomrule
\end{tabular}
\end{table}

The larger budget makes the progression especially clear.  At the
\(0.50\)- and \(0.60\)-quantiles, pure CV is the stronger of the two
constituent variance-reduction methods, and AU-CVLR remains
essentially at the CV end of the spectrum.  At the \(0.70\)- and
\(0.80\)-quantiles, modest measure adaptation begins to complement
the controls.  At the \(0.90\)-quantile, the relative ordering of the
constituent methods reverses: pure LR achieves a variance-reduction
factor of \(1.74\), compared with \(1.28\) for pure CV, while
AU-CVLR increases the factor further to \(1.93\).

At the \(0.95\)-quantile the measure component is dominant, and by
the \(0.99\)-quantile pure LR and AU-CVLR are essentially
indistinguishable, with variance-reduction factors of approximately
\(4.15\) and \(4.14\), respectively.

At the \(0.999\)-quantile the larger pilot contains roughly ten
exceedances on average rather than one.  Measure adaptation is
therefore feasible in most replications, and the LR-based procedures
again obtain a substantial improvement over CMC.  Pure LR achieves a
variance-reduction factor of approximately \(3.15\), while AU-CVLR
achieves a factor of approximately \(4.41\).  This contrast between
the two budget levels emphasizes that a potentially powerful change
of measure is useful only when the pilot itself contains enough
information to learn it.

Taken together, the eight thresholds display a continuous
progression:
\[
\text{CV-dominated}
\quad\longrightarrow\quad
\text{joint CVLR}
\quad\longrightarrow\quad
\text{LR-dominated}.
\]
The two mechanisms address different sources of inefficiency.
Control variates exploit predictable variation in the simulated
output represented by the chosen control family, whereas
likelihood-ratio adaptation changes where simulation effort is
spent.  AU-CVLR allows the variance-reduction burden to move between
these mechanisms according to the estimand and the amount of
training information available.

\begin{remark}[AU-CVLR's Implementation Complexity]
AU-CVLR's gains come at the expense of deploying additional
computational machinery in service of the LR-side update.
Unsurprisingly, this leads to an implementation that is less
straightforward than pure CV, and far less straightforward than CMC.
Implementation complexity is therefore an argument against AU-CVLR
in settings known to have strong controls, where pure CV is likely
to dominate pure LR.
\end{remark}

All reported finite-budget variances are estimated across independent
macroreplications of the complete procedure.  Simulation used solely
to assess the reported variances is not included in the method's
budget and is never fed back into fitting, tuning, or budget
allocation.  Consequently, every simulated observation that can
influence the constructed estimator is charged to that method,
whereas the independent evaluation simulation serves only to compare
the resulting procedures.

\subsection{Effect of Problem Regularization}
\label{sec:lambda-path}

We next examine the role of the regularization parameter $\lambda$ (see~\eqref{vrmaster}) on the
CVLR design.  Recall that the CVLR design kernel is
\[
\simout_\beta(x,Q)
=
R_\beta(x)\frac{dP_0}{dQ}(x),
\qquad
R_\beta(x)
=
\simout_0(x)-\beta^\top\{g(x)-\gamma\},
\]
and that the population design problem is
\[
\min_{\beta,\;Q}
\left\{
\E_Q[\simout_\beta(X,Q)^2]
+
\lambda D_{\rm KL}(Q\|P_0)
\right\}.
\]
Thus, $\lambda$ controls the extent to which the sampling distribution is
allowed to depart from the native law $P_0$.  Small values of $\lambda$
permit aggressive measure adaptation, whereas large values increasingly
penalize departures from $P_0$.

There is a simple limiting connection with classical control variates.
Let $(\beta_\lambda^\ast,Q_\lambda^\ast)$
denote an optimal CVLR design for a given $\lambda$, and let $\beta_{\rm CV}
\in
\arg\min_\beta
\E_{P_0}[R_\beta(X)^2]$
denote the optimal pure-CV coefficient.  Define $M_\lambda
=
\E_{Q_\lambda^\ast}
\left[
\simout_{\beta_\lambda^\ast}(X,Q_\lambda^\ast)^2
\right]$
and $M_{\rm CV}
=
\E_{P_0}
\left[
R_{\beta_{\rm CV}}(X)^2
\right].$
Since $(\beta_{\rm CV},P_0)$ is feasible for the CVLR problem and
$D_{\rm KL}(P_0\|P_0)=0$,
\[
M_\lambda
+
\lambda
D_{\rm KL}(Q_\lambda^\ast\|P_0)
\le
M_{\rm CV}.
\]
All designs under consideration are unbiased for $\theta$, and hence
$M_\lambda\ge\theta^2$.  Consequently,
\[
D_{\rm KL}(Q_\lambda^\ast\|P_0)
\le
\frac{M_{\rm CV}-\theta^2}{\lambda}
=
\frac{V_{\rm CV}}{\lambda},
\tag{20}
\]
where $V_{\rm CV}$ denotes the variance of the optimal pure-CV estimator.
Therefore,
\[
D_{\rm KL}(Q_\lambda^\ast\|P_0)\longrightarrow0
\qquad\text{as}\qquad
\lambda\longrightarrow\infty.
\]
In particular, Pinsker's inequality gives
\[
\|Q_\lambda^\ast-P_0\|_{\rm TV}
\le
\sqrt{\frac{V_{\rm CV}}{2\lambda}}.
\]
Under the usual continuity and uniqueness conditions, this implies $Q_\lambda^\ast\longrightarrow P_0,
\beta_\lambda^\ast\longrightarrow\beta_{\rm CV},$
and hence
\[
\mbox{Var}_{Q_\lambda^\ast}
\left[
\simout_{\beta_\lambda^\ast}(X,Q_\lambda^\ast)
\right]
\longrightarrow
V_{\rm CV}.
\tag{21}
\]
Thus, classical control variates arise naturally as the
$\lambda\to\infty$ endpoint of the CVLR regularization path.

This interpretation is also visible directly in the measure
update for CVLR.  Writing $r_k=\frac{dQ_k}{dP_0},$
the update has the form $r_{k+1}(x)
\propto
r_k(x)^{1/(1+\alpha_k\lambda)}
\exp\left\{
-\frac{\alpha_k}{1+\alpha_k\lambda}
\psi_{k+1}(x)
\right\}.$ For fixed $\alpha_k>0$,
$\frac{1}{1+\alpha_k\lambda}\longrightarrow0,
\frac{\alpha_k}{1+\alpha_k\lambda}\longrightarrow0$
as $\lambda\to\infty$, so that the measure update is increasingly pulled
toward $r_{k+1}\equiv1$, or equivalently $Q_{k+1}=P_0$.

To illustrate this behavior numerically, we repeated the $M/M/1$
queue experiment %
while varying $\lambda$ and holding all remaining
algorithmic parameters fixed.  We used
$\alpha_k\equiv0.004$ and number of iterations $K=8$.  Table~\ref{tab:lambda-path}
reports the ratio $\frac{V_{\rm CVLR}(\lambda)}{V_{\rm CV}},$
together with the estimated KL divergence between the final sampling law and
$P_0$.  A variance ratio below one favors the CVLR design, whereas convergence
of the ratio to one corresponds to recovery of the pure-CV estimator.

\begin{table}[t]
\centering
\caption{Effect of $\lambda$ on the CVLR design for the M/M/1 example.  Variance is reported
relative to the optimal pure-CV variance.}
\label{tab:lambda-path}
\begin{tabular}{rcc}
\hline
$\lambda$
& $V_{\rm CVLR}/V_{\rm CV}$
& $D_{\rm KL}(Q_K\|P_0)$
\\
\hline
0     & 0.297 & 0.0494              \\
1     & 0.303 & 0.0478              \\
10    & 0.357 & 0.0359              \\
50    & 0.539 & 0.0134              \\
100   & 0.665 & 0.00588             \\
250   & 0.814 & 0.00147             \\
500   & 0.892 & 0.000449            \\
1000  & 0.940 & 0.000128            \\
2500  & 0.974 & $2.2\times10^{-5}$  \\
5000  & 0.987 & $5.8\times10^{-6}$  \\
10000 & 0.993 & $1.5\times10^{-6}$  \\
\hline
\end{tabular}
\end{table}

The numerical behavior closely follows the population interpretation.
For small $\lambda$, the CVLR design allows substantial movement away from
$P_0$ and achieves considerably smaller variance than pure CV.  As
$\lambda$ increases, the KL divergence decreases monotonically and the
variance ratio moves steadily toward one.  In the queueing example, the
CVLR variance is approximately $30\%$ of the pure-CV variance when
$\lambda$ is near zero, but rises to approximately $99.3\%$ of the
pure-CV variance by $\lambda=10{,}000$.  %

The scale of $\lambda$ should be interpreted together with the
mirror-descent stepsize.  In the AU-CVLR mirror descent recursion, the two parameters
appear through the quantities $\frac{1}{1+\alpha_k\lambda}\text{ and }
\frac{\alpha_k}{1+\alpha_k\lambda}.$
Consequently, apparently large numerical values of $\lambda$ need not
correspond to extreme regularization when $\alpha_k$ is small.

The regularization path provides a useful interpretation of the proposed
method.  Rather than viewing pure CV and joint CVLR design as unrelated
variance-reduction procedures, they can be regarded as points on the same
continuum.  At one end, small $\lambda$ permits substantial adaptation of
both the estimator and the sampling measure.  At the other,
$\lambda\to\infty$ suppresses measure adaptation and recovers classical
control variates.  The parameter $\lambda$ therefore determines how much of
the variance-reduction burden is assigned to sampling-law adaptation beyond
what can already be achieved by the control-variate component.

\subsection{Community First Responders}

Community first responder (CFR) schemes are designed to ensure rapid response by a volunteer to the most urgent of calls to the emergency medical services (EMS), mostly out-of-hospital cardiac arrest (OHCA) \citep{berhenjagli21,beretal24b,beretal24}. Volunteers (CFRs) are alerted by a smartphone app according to a policy $\pi$ that depends on estimated travel times of all CFRs within a certain response radius of the patient. At the time the system is activated for an OHCA, one observes the estimated travel times, $t_1 < t_2 < \cdots < t_k$ for each of the $k$ CFRs within the response radius, along with the travel time, $a > t_k$ say, of the nearest available ambulance, with all times measured in minutes. The policy then decides which CFRs to alert and at what times. An alerted CFR responds to the alert with a yes/no decision after a delay that is gamma distributed with mean 1 minute and variance 1/2 minute$^2$, so it has shape parameter 2 and scale parameter 1/2 minute. CFRs need not accept an alert; we assume the probability of acceptance is $p(t)$ for a CFR with an estimated travel time $t$ from the patient, where the function $p(\cdot)$ is decreasing. For this example we take $p(t) = 0.5 e^{-t/6}$. Moreover, we assume that CFR acceptance decisions are independent of one another.

The goal is to maximize the survival of the patient, as measured through a survival function $s$ that is a function of the response time of the first CFR to arrive to the patient, whether that CFR is carrying a defibrillator, and the ambulance response time. For simplicity, however, in this example we use as the performance measure the expected time until either a CFR or the ambulance first reaches the patient. We consider a policy $\pi$ that initially alerts the three closest CFRs. Thereafter, whenever an active alert is no longer pending, the policy alerts another CFR so as to maintain three active alerts, provided no CFR has yet accepted, the ambulance has not yet arrived, and the next CFR to be alerted could still reach the patient before the ambulance if the alert were accepted immediately.

For control variates we use $C_j$, which is equal to
the indicator that CFR $j$ is the closest CFR to accept the alert, minus its mean, both scaled by the $j$th response time $t_j$. If $I_j (\bar I_j)$ denotes the indicator that the $j$th CFR accepts the alert (its complement) and $p(t_j) (\bar p(t_j))$ denotes the probability of accepting (not accepting), then
$$C_j = t_j [\bar I_1 \bar I_2 \cdots \bar I_{j-1} I_j - \bar p(t_1) \bar p(t_2) \cdots \bar p(t_{j-1}) p(t_j)].$$
This expression assumes all of the indicators $I_1, I_2, \ldots, I_k$ are generated, even if fewer than $k$ CFRs are alerted.

We take $\mathcal{Q}$ to be the family of distributions under which the acceptance decision of CFR $j$ is made with probability $q(t_j)$, independently of all other CFRs. We do not change the distribution of the time for CFRs to respond to alerts. %

For our experiment we take the ambulance travel time $a=9$ minutes, $k=8$ CFRs and the travel times of the responders to be
$$(t_1, t_2, \ldots, t_8) = (1, 3, 4, 5, 6, 6.5, 7, 7.5).$$
In practice, these values would be realized differently for every OHCA, but we fix them in our experiment. We could also have defined an experiment where the ambulance travel time, number of CFRs and their travel times were all sampled on each replication, but the present framing is simpler while still being realistic.

We again compare crude MC, pure-CV , pure-LR and AU-CVLR. The last two methods use the same measure-side update, which involves remaining within the class of independent Bernoulli distributions.

For the implementation parameters, we use
\[
K=64,\qquad
\lambda=1,
\qquad
\alpha_k\equiv0.005.
\]
For every design iteration, $n_k=4{,}000$
fresh native observations are used for the control-variate update and $m_k=8{,}000$
independent native observations are used to estimate the normalizing
constant in the measure update.  After constructing $Q_K$, we perform one
final estimator-side update using $n_{K}=4{,}000$
additional native observations.  This produces the final joint design $(\widehat\beta_{K},\widehat Q_K).$ The variances reported below are estimated using a separate evaluation
sample of $N_{\rm eval}=500{,}000$
independent observations from $P_0$. (As in the M/M/1 example, the second moment $\E_{Q_K}
\left[
\simout_\beta(X,Q_K)^2
\right]$ of the final design is evaluated
using a large independent post-experiment native sample $N_{\rm eval}$.)

For the most accurate estimator reported here, a 95\% confidence interval for the estimated mean response time was
\[
\widehat\theta=4.819 \pm 0.003\text{ minutes}.
\]
Table~\ref{tab:cfr-results} reports the estimated one-sample variance and the
corresponding variance-reduction factor relative to crude Monte Carlo.

\begin{table}[t]
\centering
\caption{Variance reduction for the CFR problem.}
\label{tab:cfr-results}
\begin{tabular}{lcc}
\hline
Method & Estimated variance & VR factor \\
\hline
Crude Monte Carlo     & 7.11 & 1.00 \\
Pure control variates & 1.24 & 5.7 \\
Pure likelihood ratio & 2.75 & 2.6 \\
CVLR Design Kernel        & 1.11 & 6.4 \\
\hline
\end{tabular}
\end{table}

Pure control variates reduce the
variance by a factor of approximately $5.7$, whereas adapting only the
sampling measure yields a factor of approximately $2.6$.  Joint adaptation
produces a modestly larger factor relative to pure control variates of approximately $6.4$. Thus, AU-CVLR yields modest gains beyond the better of control variates and importance sampling.

Similar comments on implementation and computational complexity apply in this example, as they did for the M/M/1 example.

\subsection{Option Pricing}
\label{sec:barrier-asian}

We next consider the pricing of a discretely monitored down-and-in
arithmetic Asian put under the Black--Scholes model.  This example is
designed to illustrate a setting in which the control-variate and
importance-sampling components address distinct sources of simulation
variance.  The barrier feature makes the payoff sensitive to relatively
unusual downward excursions of the asset-price path, while the arithmetic
average introduces substantial path-dependent variation even conditional
on barrier activation.

Under the native risk-neutral measure \(P_0\), the asset price satisfies
\[
\frac{dS_t}{S_t}
=
r\,dt+\sigma\,dW_t,
\qquad
S_0>0,
\]
where \(r\) is the risk-free rate, \(\sigma>0\) is the volatility, and
\(W\) is a standard Brownian motion under \(P_0\).  Equivalently,
\[
S_t
=
S_0
\exp\left\{
\left(r-\frac{\sigma^2}{2}\right)t+\sigma W_t
\right\}.
\]

Let
\[
0<t_1<\cdots<t_m=T
\]
denote the monitoring times, with \(t_j=jT/m\) in the numerical
experiments.  Define the arithmetic and geometric averages
\[
A_m
=
\frac1m\sum_{j=1}^m S_{t_j},
\qquad
G_m
=
\exp\left\{
\frac1m\sum_{j=1}^m\log S_{t_j}
\right\},
\]
and the discretely monitored running minimum
\[
M_m
=
\min_{1\le j\le m}S_{t_j}.
\]
The target quantity is the time-zero price
\[
\theta
=
\E_{P_0}[Y_0],
\]
where the discounted payoff is
\begin{equation}
Y_0
=
e^{-rT}
\mathbf 1\{M_m\le B\}
(K-A_m)^+.
\label{eq:barrier-asian-payoff}
\end{equation}
Thus, the option pays as an arithmetic Asian put only if the asset has
crossed the down barrier \(B\) at one of the monitoring times. For our experiments, we use
\[
S_0=100,
\qquad
r=0.02,
\qquad
\sigma=0.25,
\qquad
T=1,
\]
with weekly monitoring \(m=52\), strike \(K=90\), and barrier \(B=70\). 

We use three controls, all of whose native expectations are available in
closed form under the Black--Scholes model.  Define
\begin{align}
g_1
&=
e^{-rT}(K-S_T)^+,
\label{eq:barrier-asian-control-euro-put}\\
g_2
&=
e^{-rT}(K-G_m)^+,
\label{eq:barrier-asian-control-geo-asian}\\
g_3
&=
e^{-rT}S_T.
\label{eq:barrier-asian-control-stock}
\end{align}
The first control is a European put having the same strike and maturity as
the target option.  The second is a geometric-average Asian put, which is
strongly correlated with the arithmetic-average payoff but has a
closed-form price because \(\log G_m\) is Gaussian.  The final control is
the discounted terminal stock price. Writing \(g=(g_1,g_2,g_3)^\top\), the corresponding known mean vector is
\[
\gamma
=
\E_{P_0}[g]
=
\begin{pmatrix}
p_{\mathrm{Euro}}, & 
p_{\mathrm{GeoAsian}}, &
S_0
\end{pmatrix}^\trans,
\]
where $p_{\mathrm{Euro}}$ and $p_{\mathrm{GeoAsian}}$ are the standard Black-Scholes closed-form prices of a European put and a discretely monitored geometric-average Asian put.

The controls in~\eqref{eq:barrier-asian-control-euro-put}--%
\eqref{eq:barrier-asian-control-stock} are computationally inexpensive:
the terminal price, arithmetic average, and geometric average are all available once an asset-price path has been simulated.
They are intended to remove ordinary terminal-price and average-price
variation, while leaving the measure change to target the rarer
barrier-activation mechanism. We use the (standard) method described earlier for choosing the control variate coefficient vector $\beta$.

We use a two-parameter family of importance-sampling distributions that
changes the Brownian drift separately over the first and second halves of
the option lifetime. We assume, for simplicity, that the number \(m\) of monitoring times is
even.  Under the importance-sampling measure
\(Q_{\vartheta_1,\vartheta_2}\), the Brownian increments satisfy
\[
dW_t
=
\begin{cases}
\vartheta_1\,dt+d\widetilde W_t,
&
0\le t\le T/2,
\\[4pt]
\vartheta_2\,dt+d\widetilde W_t,
&
T/2<t\le T,
\end{cases}
\]
where \(\widetilde W\) is a standard Brownian motion under
\(Q_{\vartheta_1,\vartheta_2}\).  Equivalently, the simulated asset path
under \(Q_{\vartheta_1,\vartheta_2}\) is
\[
S_t
=
S_0
\exp\left\{
\left(r-\frac{\sigma^2}{2}\right)t
+
\sigma W_t
\right\},
\]
where \(W_t\) is generated with the piecewise constant drifts specified
above.

Let $W^{(1)} = W_{T/2}-W_0$ and $W^{(2)}=W_T-W_{T/2}$
denote the Brownian increments over the first and second blocks under the
sampling measure.  The likelihood ratio is
\begin{align}
L_{\vartheta_1,\vartheta_2}
&=
\frac{dP_0}{dQ_{\vartheta_1,\vartheta_2}}
\nonumber\\
&=
\exp\left\{
-\vartheta_1W^{(1)}
-\vartheta_2W^{(2)}
+
\frac12\vartheta_1^2\frac T 2
+
\frac12\vartheta_2^2 \frac T 2
\right\}.
\label{eq:two-block-lr}
\end{align}

The importance-sampling estimator is therefore
\[
T_{\rm IS}
=
Y_0 L_{\widehat\vartheta_{1,\rm IS},
       \widehat\vartheta_{2,\rm IS}},
\qquad
X\sim
Q_{\widehat\vartheta_{1,\rm IS},
  \widehat\vartheta_{2,\rm IS}}.
\]

Negative values of \(\vartheta_1\) encourage an early downward movement in
the asset price and hence increase the frequency of barrier activation.
Negative values of \(\vartheta_2\) encourage low values of the arithmetic
average and terminal stock price conditional on the early path history.

We use an independent native pilot sample
\[
\{(Y_i,g_i,W_i^{(1)},W_i^{(2)}):i=1,\ldots,n_{\rm pilot} = 10,000\}
\]
generated under \(P_0\) that is split into a fitting
sample (60\% of the pilot sample) and a validation sample (40\%).

Our choice of the importance sampling family reduces the need for regularization, so we do not regularize. We could use a gradient search to choose the two importance sampling parameters, but due to the low dimension of the parameter space we instead simply consider a grid of drift pairs, where the two drifts take values in the sets
\begin{align*}
\mathcal V_1 &=
\{-4.0,-3.5,-3.0,-2.5,-2.0,-1.5,-1.0,-0.5,0\},\\
\mathcal V_2 &=
\{-3.0,-2.5,-2.0,-1.5,-1.0,-0.5,0\}.
\end{align*}
The larger range of negative first-block drift values reflects the greater role of the first block in generating barrier crossings.

For pure importance sampling, we set \(\beta=0\) and select $(\vartheta_1, \vartheta_2)$ to be the pair from the grid that minimizes the estimated
second moment
\[
\E_{Q_\vartheta}[Y_0^2L_\vartheta^2]
=
\E_{P_0}[Y_0^2L_\vartheta].
\]

For the joint CVLR design, we select $\vartheta_1, \vartheta_2$ to minimize the estimated second moment of the control-variate estimator using the $n_{\rm val}$ validation pilot observations
\begin{align}
\widehat\vartheta_{\rm CVLR}
&\in
\arg\min_{\vartheta\in\mathcal V}
\left\{
\frac{1}{n_{\rm val}}
\sum_{i=1}^{n_{\rm val}}
\left[
Y_i-
\widehat\beta_\vartheta^\top(g_i-\gamma)
\right]^2
L_{\vartheta,i}
\right\},
\label{eq:option-cvlr-grid-search}
\end{align}
where $L_\vartheta,i$ is the likelihood ratio of the $i$th sample.
After selecting $\vartheta_1, \vartheta_2$ in this manner,
the corresponding control coefficient is refitted using all pilot
observations.
Likewise, the pure-CV coefficient is estimated from the complete pilot.

The production samples used to evaluate the four estimators are independent
of the pilot sample.  Thus, pilot observations influence estimator design
but are not reused in the final production averages.

For each method, let \(Y\) denote its one-sample estimator.  We report its
estimated one-sample variance, $\widehat V = \widehat{\operatorname{Var}}(Y)$,
and the usual variance-reduction factor relative to crude Monte Carlo, $\widehat{\mathrm{VRF}}
= \widehat V_{\rm CMC} / \widehat V$.

To account for method-dependent computational effort, let
\[
\widehat c_{\rm prod}
=
\frac{\text{production wall-clock time}}
{\text{number of retained production observations}}
\]
denote the measured production time per observation and define the
production-cost-adjusted variance as
\begin{equation}
\widehat V_{\rm prod}
=
\widehat c_{\rm prod}\widehat V.
\label{eq:option-production-adjusted-variance}
\end{equation}
Smaller values of~\eqref{eq:option-production-adjusted-variance}
correspond to greater asymptotic computational efficiency.

We also report an amortized cost-adjusted variance that charges the
pilot-design effort to the procedure.  Let \(C_{\rm design}\) denote the
wall-clock time used to generate the pilot and select the estimator design,
and let \(n_{\rm prod} = 500,000\) denote the number of production observations.
Then
\[
\widehat c_{\rm amort}
=
\widehat c_{\rm prod}
+
\frac{C_{\rm design}}{n_{\rm prod}},
\]
and the amortized cost-adjusted variance is
\[
\widehat V_{\rm amort}
=
\widehat c_{\rm amort}\widehat V.
\label{eq:option-amortized-adjusted-variance}
\]
This latter measure is conservative for applications in which a learned
importance-sampling design is reused repeatedly, but it provides a
transparent accounting of the computational cost of adaptation for a
single pricing calculation.

The corresponding relative production and amortized efficiencies are
\[
\widehat{\mathrm{Eff}}_{\rm prod}
=
\frac{
\widehat V_{\rm prod,CMC}
}{
\widehat V_{\rm prod}
},
\qquad
\widehat{\mathrm{Eff}}_{\rm amort}
=
\frac{
\widehat V_{\rm amort,CMC}
}{
\widehat V_{\rm amort}
}.
\]
Values greater than one favor the corresponding method over crude Monte
Carlo.

Table~\ref{tab:barrier-asian-results} summarizes the results.  The price estimates themselves all agree to statistical error, with a 95\% confidence interval of $1.285 \pm 0.003$ for the most accurate estimator. However, the primary comparison concerns the one-sample
variance and the two cost-adjusted variance measures.

\begin{table}[t]
\centering
\small
\caption{Variance and computational-efficiency comparison for the
discretely monitored down-and-in arithmetic Asian put.}
\label{tab:barrier-asian-results}
\begin{tabular}{lrrrrr}
\toprule
Method
& Variance
& Factor
& $\widehat c_{\rm prod}$
& $\widehat{\mathrm{Eff}}_{\rm prod}$
& $\widehat{\mathrm{Eff}}_{\rm amort}$
\\
& $\widehat{\mathrm{V}}$ &$\widehat{\mathrm{VRF}}$ &$10^{-6}$ sec./obs& &
\\
\midrule
CMC
& \texttt{14.6}
& \(1\)
& \texttt{1.11}
& \(1\)
& \(1\)
\\
Pure-CV
& \texttt{1.53}
& \texttt{9.6}
& \texttt{1.15}
& \texttt{9.3}
& \texttt{8.9}
\\
Pure-LR
& \texttt{1.97}
& \texttt{7.4}
& \texttt{1.15}
& \texttt{7.2}
& \texttt{6.7}
\\
CVLR
& \texttt{1.00}
& \texttt{14.6}
& \texttt{1.14}
& \texttt{14.3}
& \texttt{13.3}
\\
\bottomrule
\end{tabular}
\end{table}

The controls remove variation
associated with broad terminal-price and average-price movements, but
they do not completely account for the path-dependent discrepancy between
a barrier-triggered arithmetic Asian payoff and its smooth proxy payoffs.
Importance sampling then reallocates probability mass
toward paths with an early downward move and, when beneficial, continued
downward pressure during the second half of the contract life.  The joint
CVLR design can therefore obtain variance reduction beyond that delivered
by either controls or importance sampling alone.

In the present experiment, the estimated relative amortized efficiencies are 8.9, 6.7 and 13.3 for CVs, IS and CVLR respectively, indicating strong additional gains of CVLR over both CVs and IS. These values are slightly smaller than their respective relative production efficiencies, reflecting the overhead in the pilot.

\section*{Acknowledgements}
We thank Barry Nelson for many years of stimulating conversations, generous mentoring, friendship, and laughter. We are equally grateful for the clarity of his exposition, the depth of his thinking, and the creativity of his ideas for solving simulation problems.

\par\medskip
\noindent\textbf{Funding:} This research was supported by the National Science Foundation grants OAC 2410950 and DMS 2230023.
\bibliographystyle{plainnat}
\bibliography{FromShane}
\end{document}